\documentclass[a4paper,12pt]{article}
\usepackage{amscd,amssymb,amsmath,amsthm}
\usepackage{rotating}
\usepackage{hyperref}
\usepackage{graphicx,caption,dsfont}
\usepackage{mathtools}

\usepackage{subcaption}
\usepackage{amsfonts,enumerate}
\usepackage{fullpage}
\usepackage[numbers,sort & compress]{natbib}
\usepackage{soul}
\usepackage{xfrac}
\usepackage{booktabs}
\numberwithin{equation}{section}

\newtheorem{thm}{Theorem}[section]

\newtheoremstyle{case}{}{}{}{}{}{:}{ }{}
\theoremstyle{case}

\allowdisplaybreaks

\makeatletter
\def\and{%
  \end{tabular}%
  \begin{tabular}[t]{c}}%
\def\@fnsymbol#1{\ensuremath{\ifcase#1\or 1\or 2\or 3\or
   4\or 5\or 6\or 7\or 8\or 9\else\@ctrerr\fi}}
\makeatother
\vspace{-8ex}
\date{}
\begin{document}
\baselineskip=0.71cm
\begin{center}
{\large \textbf{Initial value problem for the two-dimensional time-fractional generalized convection-reaction-diffusion-wave equation: Invariant subspaces and exact solutions}}
\end{center}
\vspace{-0.5cm}
\begin{center}
P. Prakash$^{*}$,  K.S. Priyendhu,  K.M. Anjitha\\
Department of Mathematics,\\
Amrita School of Engineering, Coimbatore-641112,\\
  Amrita Vishwa Vidyapeetham, INDIA.\\
E-mail ID:vishnuindia89@gmail.com\\
\hspace{1.6cm} :p$_{-}$prakash@cb.amrita.edu\\
$*$-Corresponding Author
\end{center}
\begin{abstract}
This work investigates how we can extend the invariant subspace method to two-dimensional time-fractional non-linear PDEs.
More precisely, the systematic study has been provided for constructing the various dimensions of the invariant subspaces for the two-dimensional time-fractional generalized convection-reaction-diffusion-wave equation along with the initial conditions for the first time. Additionally, the special types of the above-mentioned equation are discussed through this method separately such as reaction-diffusion-wave equation, convection-diffusion-wave equation and diffusion-wave equation. Moreover, we explain how to derive variety of exact solutions for the underlying equation along with initial conditions using the obtained invariant subspaces. Finally, we extend this method to two-dimensional time-fractional non-linear PDEs with time delay.  Also, the effectiveness and applicability of the method have been illustrated through the two-dimensional time-fractional cubic non-linear convection-reaction-diffusion-wave equation with time delay. In addition, we observe that the obtained exact solutions can be viewed as the combinations of Mittag-Leffler function and polynomial, exponential and trigonometric type functions.
\end{abstract}
\section{Introduction}
 In recent years, the study of fractional differential equations (FDEs) has gained considerable popularity and importance due to the exact description of widespread applications in various areas of science and engineering \cite{pi,kai,kts,vas5,mo11,mac} such as physics \cite{vas5}, mechanics \cite{mn}, electrodynamics \cite{vas7}, biology \cite{vas14}, visco-elasticity \cite{btr} and so on \cite{vas2}. Unlike the integer-order differential operator, fractional-order (non-integer order) differential operators are non-local and it depends not only the immediate past but also on all the past values of the function. Thus, the fractional derivative operators are typically used to add memory in a complex system. Hence the fractional-order models help to understand the behaviour of the system and phenomena that are classified by power-law non-locality, power-law long term memory and fractal properties \cite{vas2}. Also, we know that many complex systems depend not only on an instantaneous time but also on all the past time history of the system which can be successfully modelled by using the theory of non-integer order derivatives and non-integer order integrals \cite{pi,kai,kts,vas5,mo11,mn,vas7,vas14,btr,vas2,mac}. In addition, it is important to note that one of the advantages of fractional PDEs is help to study both the hyperbolic and parabolic types of PDEs simultaneously. For example, let us consider the fractional-order PDE
\begin{equation}\label{d}
\dfrac{\partial^\alpha u}{\partial t^\alpha}=c\dfrac{\partial^2u}{\partial x^2},\ c\in\mathbb{R},\ t>0,\ \alpha\in(0,2].
\end{equation}
Note that the above equation \eqref{d} is referred to as a time-fractional diffusion-wave equation \cite{bb1,lk,wu2}. Additionally, we observe that
\begin{itemize}
\item[(i)] If $\alpha=1$, then the equation \eqref{d} is popularly known as the standard diffusion equation.
\item[(ii)] When $\alpha=2$, equation \eqref{d} is familiarly called a standard wave equation.
\item[(iii)] If $\alpha<1$, equation \eqref{d} describes the sub-diffusion phenomenon \cite{bb1,lk,wu} which represents the slow movement of particles. Also, if $\alpha>1$, equation \eqref{d} can be viewed as a super-diffusion phenomenon \cite{bb1} that represents the fast movement of particles.
\item[(iv)] It is important to note that when $\alpha\in(1,2)$, equation \eqref{d} helps to study the intermediate process between diffusion and wave phenomena \cite{bb1}. For this case, the above equation \eqref{d} is known as fractional diffusion-wave equation \cite{bb1,lk}.
\end{itemize}

It is well-known that the generalized non-linear convection-reaction-diffusion (CRD) equations \cite{cb,CRD1,CRD2,CRD3,CRD4,CRD5,ps1} have a wide range of applications in the fields of biology, physics, chemistry and engineering such as heat transfer processes, adsorption in porous medium, chemical reactions, population dynamics, lubrication and viscosity of fluids, turbulence of heterogeneous physical system, amount of contamination transported through ground water and so on. The CRD equation can be classified into convection, reaction and diffusion processes. Convection refers to the movement of a substance from one region to another, reaction process is due to adsorption, decay and reaction of substances with other components. Finally, the key process is diffusion which is due to the movement of substances from an area of high concentration to an area of low concentration. The CRD equation is a mathematical model that describe how the concentration of the substance distributed in the medium changes under the influence of the above-mentioned three processes. Also, we note that the generalized CRD equation can be categorized into three important special types that are (i) diffusion, (ii) reaction-diffusion and (iii) convection-diffusion or absorption-diffusion phenomena.
For instance, the Newell-Whitehead equation or amplitude equation \cite{CRD4}, the Fisher-KPP equation or the logistic equation \cite{CRD4}, the FitzHugh-Nagumo equation or the bistable equation \cite{CRD4} and the Zeldovich equation or Huxley equation \cite{CRD4} are reaction-diffusion equations. The Richards equation \cite{CRD4}, the foam drainage equation \cite{CRD4} and the Burgers equation \cite{CRD4} are non-linear convection-diffusion equations.

We know that fractional-order derivative operators violate some standard fundamental properties such as the semigroup property, chain rule and Leibniz rule. The violation of the above-mentioned properties is the fundamental properties of the fractional-order derivatives that allow us to describe the memory \cite{vas2}.
Using the fractional-order derivative, the present article deals with the systematic study for finding the exact solutions of the initial value problem for the two-dimensional time-fractional generalized non-linear convection-reaction-diffusion-wave equation in the form
	\begin{equation}\label{1}
	\frac{\partial^{\alpha}u}{\partial t^{\alpha}}=
 \sum\limits^{2}_{i=1} \frac{\partial }{\partial x_i}\left(A_i(u)\frac{\partial u}{\partial x_i}\right)+ \sum\limits^{2}_{i=1}B_i(u)\frac{\partial u}{\partial x_i}+C(u),\ \alpha\in(0,2],
 \end{equation}
along with the initial conditions
\begin{itemize}
\item[$\bullet$] $u(x_1,x_2,0)= \varphi_1(x_1,x_2)$  if $\alpha\in(0,1]$,
\item[$\bullet$] $u(x_1,x_2,0)=\varphi_1(x_1,x_2)$ and ${\dfrac{\partial u}{\partial t}}\mid_{t=0}=\varphi_2(x_1,x_2)$ if $\alpha \in(1,2]$,
\end{itemize}
 where $u=u(x_1,x_2,t)$ and $\dfrac{\partial^{\alpha}u}{\partial t^{\alpha}}$ denotes the Caputo fractional derivative of order $\alpha>0$, is defined by \cite{kai,pi}
 \begin{eqnarray}
\begin{aligned}\label{c}
\dfrac{\partial^\alpha u}{\partial t^\alpha}=
\left\{
                                            \begin{array}{ll}
                                              \dfrac{1}{\Gamma(m-\alpha)}\int_{0}^{t}\dfrac{\partial ^{m} u(x_1,x_2,\tau)}{\partial \tau^{m} }(t-\tau)^{m-\alpha-1}d\tau, & m-1<\alpha<m,\ m\in\mathbb{N}; \\
                                              \dfrac{\partial^m u}{\partial t^m}, & \alpha=m\in\mathbb{N}
                                            \end{array}
                                          \right.
\end{aligned}
\end{eqnarray}
  and the functions $ A_i(u)$ denote the process of  diffusion, $ B_i(u)$ denote convection, and $C(u)$ denotes reaction of a system under consideration with space variables $x_1, x_2$ and time-variable $t$.

In the literature, there is no well-defined method for solving the FDEs due to the violations of the fundamental properties of fractional-order derivatives. Hence finding the exact solutions of fractional non-linear differential equations are a very challenging and toughest job. However, the derivation of exact solutions of FDEs is an important task because it will help to understand the behaviour of the complex systems. Due to these reasons, in recent years, many mathematicians and physicists have spent much attention on developing the exact and numerical methods for solving the FDEs such as invariant subspace method \cite{ra,hes,pra4,pra2,pra3,ps2,ru2,pa1,s1,pp1,pp2,ru3,praf}, Lie group method \cite{praf,pra2,ru4,pra1,kdv,ba,pra5,nas,ppp,ms2}, element free Galerkin (EFG) method \cite{galer}, quadratic spline collocation method \cite{wu}, Adomian decomposition method \cite{wu1,Ge,mo}, differential transform method \cite{mo2}, iterative method \cite{abd1}, new hybrid method \cite{mma}, difference methods \cite{abd}, variational approach \cite{wu2} and so on. Recent studies have shown that the invariant subspace method is a very effective and powerful analytical approach for deriving the exact solutions of fractional non-linear PDEs.

The invariant subspace method is also called the generalized separation of variables method. This method was initially studied by Galaktionov and Svirshchevskii \cite{gs}.  After that, it was developed for scalar and coupled non-linear PDEs by Ma \cite{wx} and many others \cite{wx1,wy,yem,li}. Recently, the invariant subspace method was extended for scalar and coupled non-linear fractional PDEs \cite{praf,ra,pra4,pra2,pra3,ps2,ru2,pa1,s1,pp1,pp2,ru3,hes}. Zhu and Qu \cite{zq} have extended this method to two-dimensional non-linear PDEs. However, to the best of our knowledge, no one has been extended the invariant subspace method to two-dimensional fractional non-linear PDEs. Hence the main purpose of the article is to show how to extend the invariant subspace method for two-dimensional time-fractional non-linear PDEs. More specifically, this work provides a systematic study for constructing various dimensions of the invariant subspaces for the two-dimensional time-fractional generalized convection-reaction-diffusion-wave equations with the different types of power-law non-linearities along with the initial conditions for the first time. Also, we show how to derive exact solutions for the initial value problems of the underlying equation using the obtained invariant subspaces. Additionally, this article investigates for finding the invariant subspaces for the special types of the given equation that are the convection-diffusion-wave equation, reaction-diffusion-wave equation and diffusion-wave equation.

The structure of the paper is organized as follows. In section 2, we give a detailed study for finding the invariant subspaces of the two-dimensional time-fractional non-linear PDEs. Also, we explain how to construct exact solutions for the given PDEs using the obtained invariant subspaces. In section 3, we give a brief systematic investigation for finding the two types of invariant subspaces with various dimensions for two-dimensional time-fractional convection-reaction-diffusion-wave equation with different non-linearities. In addition, the special types of the above-mentioned equation are discussed through this method separately such as reaction-diffusion-wave equation, convection-diffusion-wave equation and diffusion-wave equation. Also, we explain how to construct variety of exact solutions for the given equation along with the initial conditions using the obtained invariant subspaces in section 4.  In section 5, we explain the extension of the invariant subspaces to two-dimensional time-fractional PDEs with time delay. The applicability and effectiveness of the method have been illustrated through the two-dimensional time-fractional cubic non-linear convection-reaction-diffusion-wave equation with linear time delay. Finally, a brief summary of results and concluding remarks are given in section 6.
%
\section{Invariant subspaces admitted by two-dimensional time-fractional non-linear PDEs}
In this section, we present a detailed systematic approach for finding the invariant subspaces of two-dimensional time-fractional non-linear PDEs. Also, we give a detailed study for deriving the exact solutions of the two-dimensional time-fractional non-linear PDEs using the obtained invariant subspaces. Thus, we consider the generalized two-dimensional time-fractional non-linear PDE in the form
\begin{equation} \label{1.1}
\frac{\partial^{\alpha}u}{\partial t^{\alpha}}=\mathcal{\hat{K}}[u]
\equiv\mathcal{\hat{K}}\left( x_{1},x_{2},
\frac{\partial u}{\partial x_1},\frac{\partial u}{\partial x_2},\frac{\partial^2 u}{\partial x_1^2},\frac{\partial^2 u}{\partial x_2^2},\frac{\partial^{2}u}{\partial x_1\partial x_2},\ldots,\frac{\partial^k u}{\partial x_1^k},\frac{\partial^k u}{\partial x_2^k},\frac{\partial^{k}u}{\partial x_1^{k_1}\partial x_2^{k_2}}\right),
\end{equation}
where $u=u(x_1,x_2,t)$, $\alpha>0$, $t>0$, $x_1,x_2\in\mathbb{R}$ and $\dfrac{\partial^{\alpha}}{\partial t^{\alpha}}(\cdot) $ denotes Caputo fractional derivative \eqref{c} of order $\alpha$, and $\mathcal{\hat{K}}[u]$ is the sufficiently given smooth differential operator of order $k$, and $k_1+k_2=k$, $k_1,k_2 \in \mathbb{N}$.

Now, we assume the linearly independent set
\begin{equation}\label{1.2}
S=\{\xi_m(x_1,x_2): m=1,2,\dots,n \}.
\end{equation}
Then, we can define the n-dimensional linear space
\begin{eqnarray} \label{1.3}
\begin{aligned}
\mathcal{V}_n&=Span(S)\\
&=Span\{\xi_m(x_1,x_2): m =1,2,\dots,n\}\\
&=\left\{\sum_{m=1}^{n}\kappa_m\xi_m(x_1,x_2):\kappa_m \in \mathbb{R}, m = 1,2,\dots,n\right\},
\end{aligned}
\end{eqnarray}
where $n$ denotes the dimension of the linear space and it is denoted by $dim(\mathcal{V}_n)=n\ (<\infty)$.
The linear space $\mathcal{V}_n $ given in \eqref{1.3} is said to be invariant under the differential operator $\mathcal{\hat{K}}[u]$ given in \eqref{1.1} if for every $ u\in \mathcal{V}_n $ implies $ \mathcal{\hat{K}}[u] \in \mathcal{V}_n $, which can be written as
 \begin{equation} \label{1.4}
 \mathcal{\hat{K}}\left[\sum_{m=1}^{n}\kappa_m\xi_m(x_1,x_2) \right]=\sum \limits^{n} _{m=1}\Psi_m(\kappa_1,\kappa_2,\ldots,\kappa_n)\xi_m(x_1,x_2),
\end{equation}
where $\kappa_m \in \mathbb{R}$ and $\Psi_m $ denotes coefficient of expansion $ \mathcal{\hat{K}}[u] $ with respect to the basis set given in \eqref{1.2}, $m=1,2,\dots,n$.
\begin{thm}
Suppose that the differential operator $ \mathcal{\hat{K}}[u]$ given in \eqref{1.1} admits an n-dimensional linear space $\mathcal{V}_n$ given in \eqref{1.3},
then the underlying two-dimensional time-fractional PDE has admits an exact solution of the form
\begin{equation}\label{1.5}
u(x_1,x_2,t)=\sum_{m=1}^{n}\Phi_m(t)\xi_m(x_1,x_2),
\end{equation}
 where the functions $\Phi_m(t), m=1,2,\dots,n$ are the solutions of the following system of ODEs of fractional-order
 \begin{equation}\label{1.6}
 \frac{d^{\alpha}\Phi_m(t)}{dt^\alpha} =\Psi_m(\Phi_1(t),\Phi_2(t),...,\Phi_n(t)), \  m=1,2,\dots,n.
  \end{equation}
\end{thm}
\textbf{Proof.}
Let $\mathcal{V}_n$ be an $n$-dimensional invariant subspace admitted by the given differential operator ${\mathcal{\hat{K}}}[u(x_1,x_2,t)]$. Also, we assume that
\begin{equation}
 u(x_1,x_2,t)=\sum\limits_{m=1}^{n}\Phi_m(t)\xi_m(x_1,x_2).
 \end{equation}
 Computing the Caputo fractional derivative of $u(x_1,x_2,t)$ of order $\alpha>0$ with respect to $t$ gives
 \begin{equation}\label{1.7}
 \frac{\partial^{\alpha}u}{\partial t^{\alpha}}=\sum^{n} _{m=1}  \frac{d^{\alpha}\Phi_m(t)}{dt^\alpha}\xi_m(x_1,x_2).
 \end{equation}
 Since the linear space $\mathcal{V}_n$ is invariant under ${\mathcal{\hat{K}}}[u]$. Thus, we have
\begin{equation} \label{1.81}
\mathcal{\hat{K}}[u(x_1,x_2,t)] =\sum \limits^{n} _{m=1}\Psi_m(\Phi_1(t),\Phi_2(t),...,\Phi_n(t))\xi_m(x_1,x_2).
\end{equation}
Substituting \eqref{1.7} and \eqref{1.81} in \eqref{1.1}, we get
\begin{equation}\label{1.91}
\sum\limits_{m=1}^{n}\left[\frac{d^{\alpha}\Phi_m(t)}{dt^\alpha}-\Psi_m(\Phi_1(t),\Phi_2(t),...,\Phi_n(t))\right]\xi_m(x_1,x_2)=0.
\end{equation}
Since the functions $\xi_m(x_1,x_2)$, $i=1,2,\dots,n$ are linearly independent. Thus, the above equation \eqref{1.91} reduces to system of fractional ODEs \eqref{1.6}.
\subsection{Estimation of invariant subspaces for the given equation \eqref{1.1}}
In this subsection, we would like to explain how to find the invariant subspaces for the given differential operator $\mathcal{\hat{K}}[u]$. Thus, let us first define the following two linear spaces
\begin{equation}\begin{split}
 \mathcal{V}_{n_1} &=  Span\left\{v _1(x_1),v _2(x_1),\dots,v_{n_1}(x_1)\right\} \\
 &=\left\{y_1\ \big{|}\ \mathcal{D}_{x_1}^{n_1}[y_1]\equiv \frac{d^{n_1}y_1}{dx_1^{n_1}} +a_{n_1-1}\frac{d^{n_1-1}y_1}{dx_1^{n_1-1}}+a_{n_1-2}\frac{d^{n_1-2}y_1}{dx_1^{n_1-2}}+\dots+a_{0}y_1=0\right\},\\ \label{1.8}
 \mathcal{V}_{n_2}&=Span\left\{\zeta_1(x_2),\zeta_2(x_2),\dots,\zeta_{n_2}(x_2)\right\} \\
 &= \left\{y_2\ \big{|}\ \mathcal{D}_{x_2}^{n_2}[y_2]\equiv \frac{d^{n_2}y_2}{dx_2^{n_2}} +b_{n_2-1}\frac{d^{n_2-1}y_2}{dx_2^{n_2-1}}+b_{n_2-2}\frac{d^{n_2-2}y_2}{dx_2^{n_2-2}}+\dots+b_{0}y_2=0\right\},\\
\end{split}
\end{equation}
where $ \{v _i(x_1): i=1,2,\dots,n_1\}$ and $\{\zeta_j(x_2):j=1,2,\dots,n_2\} $ are the set of $n_k$-linearly independent solutions of the linear homogeneous ODEs $\mathcal{D}^{n_k}_{x_k}[y_k]=0$, for $ k=1,2$, respectively and $a_i,b_i\in\mathbb{R}$, $i=0,1,2,\dots,n_{i}-1,$ $i=1,2$. Here $\mathcal{D}^{n_k}_{x_k}[y_k]=0$ denotes the $n_k$-th order homogeneous linear ODEs. Now, we can formally define the type I and type II invariant linear spaces for the given non-linear differential operator $\mathcal{\hat{K}}[u]$ given in \eqref{1.1} that are discussed below in detail.
\begin{itemize}
\item[1.]\textbf{\textit{Type I linear space}}:
 First, we define the type I invariant subspace for the differential operator $\mathcal{\hat{K}}[u]$ in the form
\begin{eqnarray}
\begin{aligned}\label{1.9}
\mathcal{V}_{n_1n_2}=&\text{Span}\left\{v_1(x_1)\zeta_1(x_2),\dots,v_1(x_1)\zeta_{n_2}(x_2),\dots,v_{n_1}(x_1)\zeta_1(x_2),\dots,\right.\\
&\left.v_{n_1}(x_1)\zeta_{n_2}(x_2)\right\} \\
=&\left\{\xi(x_1,x_2)=\sum\limits_{j=1}^{n_2}\left(\sum\limits_{i=1}^{n_1}\kappa_{ij}v_i(x_1)\right)\zeta_j(x_2)\ \big{|}\ \mathcal{D}_{x_1}^{n_1}[\xi]=0, \mathcal{D}_{x_2}^{n_2}[\xi]=0\right\},
\end{aligned}
\end{eqnarray}
where $\kappa_{ij} \in \mathbb{R},\ i=1,\dots,n_1,\ j=1,2,\dots,n_2$.\\
If the differential operator $\mathcal{\hat{K}}[u]$ given in \eqref{1.1} admits the linear space $\mathcal{V}_{n_1n_2}$, then we have
 \[ \mathcal{\hat{K}}\left[\sum\limits_{i=1}^{n_1}\sum\limits_{j=1}^{n_2}\kappa_{ij}v_i(x_1)\zeta_j(x_2) \right]=\sum\limits_{j=1}^{n_2}\sum\limits_{i=1}^{n_1}\Psi_{ij}(\kappa_{11},\dots, \kappa_{n_1n_2})v_i(x_1)\zeta_j(x_2). \]
 This suggests an exact solution of the given equation \eqref{1.1} in the form
\[u(x_1,x_2,t)= \sum\limits_{j=1}^{n_2}\left(\sum\limits_{i=1}^{n_1}\Phi_{ij}(t)v_i(x_1)\right)\zeta_j(x_2) =\sum\limits_{i=1}^{n_1}\left(\sum\limits_{j=1}^{n_2}\Phi_{ij}(t)\zeta_j(x_2)\right)v_i(x_1).\]
For this case, we obtain the invariance conditions for type I linear space of the given differential operator $\mathcal{\hat{K}}[u]$ of the form
\begin{eqnarray*}\label{1.10}
&&\mathcal{D}_{x_{1}}^{n_1}[\mathcal{\hat{K}}]\equiv \frac{d^{n_1}\mathcal{\hat{K}}}{dx_1^{n_1}}+a_{n_1-1}\frac{d^{n_1-1}\mathcal{\hat{K}}}{dx_1^{n_1-1}}+\dots+a_1\frac{d\mathcal{\hat{K}}}{dx_1}+a_{0}\mathcal{\hat{K}}=0,\\
&&\mathcal{D}_{x_2}^{n_2}[\mathcal{\hat{K}}]\equiv \frac{d^{n_2}\mathcal{\hat{K}}}{dx_2^{n_2}}+b_{n_2-1}\frac{d^{n_2-1}\mathcal{\hat{K}}}{dx_2^{n_2-1}}+\dots+b_1\frac{d\mathcal{\hat{K}}}{dx_2}+b_{0}\mathcal{\hat{K}}=0
\end{eqnarray*}
along with
\begin{eqnarray*}
\begin{aligned}
&\mathcal{D}_{x_1}^{n_1}[u]\equiv \frac{d^{n_1}u}{dx_1^{n_1}}+a_{n_1-1}\frac{d^{n_1-1}u}{dx_1^{n_1-1}}+\dots+a_1\frac{d u}{dx_1}+a_{0}u=0,\\
&\mathcal{D}_{x_2}^{n_2}[u]\equiv \frac{d^{n_2}u}{dx_2^{n_2}}+b_{n_2-1}\frac{d^{n_2-1}u}{dx_2^{n_2-1}}+\dots+b_1\frac{d u}{dx_2}+b_{0}u=0
\end{aligned}
\end{eqnarray*}
and their differential consequences with respect to $x_1$ and $x_2$ if the differential operator $\mathcal{\hat{K}}[u] $ admits the type I linear space $\mathcal{V}_{n_1n_2}$.
Some important consequences of type I invariant subspaces for the given differential operator $\mathcal{\hat{K}}[u] $ are as follows
 \begin{itemize}
 \item[(i)] If $ 1 \in \mathcal{V}_{n_1} $, then $a_0 =0$ in $\mathcal{D}_{x_1}^{n_1}[\xi]=0 $.
 Suppose $ v_1(x_1)=1 \in \mathcal{V}_{n_1}$. For this case, we obtain  $\mathcal{V}_{n_1n_2}=\text{Span}\{\zeta_1(x_2),\ldots,\zeta_{n_2}(x_2),v _2(x_1)\zeta_1(x_2),\ldots,v _{n_1}(x_1)\zeta_{n_2}(x_2)\} $. Thus, we have $\text{dim}(\mathcal{V}_{n_1n_2})= n_1n_2$.
 \item[(ii)]Similarly,  $b_0 =0$  in $\mathcal{D}_{x_2}^{n_2}[\xi]=0$ when $ 1 \in \mathcal{V}_{n_2}$. Suppose $ \zeta_1(x_2)=1 \in \mathcal{V}_{n_2}$. Thus, we have $\mathcal{V}_{n_1n_2} = \text{Span}\{v _1(x_1),\ldots,v _{n_1}(x_1),v _1(x_1)\zeta_2(x_2),\ldots,v _{n_1}(x_1)\zeta_{n_2}(x_2)\} $ and  $\text{dim}( \mathcal{V}_{n_1n_2}) = n_1n_2$.
\item[(iii)] If $ 1\in \mathcal{V}_{n_1}\ \cap\ \mathcal{V}_{n_2} $, then $  a_0 =0 $ in $\mathcal{D}_{x_1}^{n_1}[\xi]=0$  and  $b_0 =0$  in $\mathcal{D}_{x_2}^{n_2}[\xi]=0$. In this case, we obtain $\mathcal{V}_{n_1n_2}=\text{Span}\{1,v _2(x_1),\dots,v _{n_1}(x_1),\zeta_2(x_2),\dots,\zeta_{n_2}(x_2),v _2(x_1)\zeta_2(x_2),$ $\dots,v _{n_1}(x_1)\zeta_{n_2}(x_2)\} $ and $\text{dim}( \mathcal{V}_{n_1n_2}) = n_1n_2.$
 \end{itemize}
\item[2.]\textbf{\textit{Type II linear space}}:
The type II invariant subspace $\mathcal{V}_{n_1+n_2-1}$ for the given differential operator $\mathcal{\hat{K}}[u]$ is defined as follows
\begin{eqnarray}\begin{aligned}
\mathcal{V}_{n_1+n_2-1}=&\text{Span}\left\{1,v_1(x_1),\dots,v_{n_1-1}(x_1),\zeta_1(x_2),\dots,\zeta_{n_2-1}(x_2)\right\} \\
\label{1.10}
=&\left\{\xi(x_1,x_2)=\kappa+\sum\limits_{i=1}^{n_1-1}\kappa_iv_i(x_1)+\sum\limits_{j=1}^{n_2-1}\delta_j\zeta_j(x_2)\ \big{|}\ \mathcal{D}^{n_1}_{x_1}[\xi]=0, \right.\\
&\left.\mathcal{D}^{n_2}_{x_2}[\xi]=0,\frac{\partial^2\xi}{\partial x_1\partial x_2}=0, a_0=b_0=0\right\}.
 \end{aligned}
\end{eqnarray}
If the given differential operator $\mathcal{\hat{K}}[u] $ admits the linear space $\mathcal{V}_{n_1+n_2-1}$, then we can write
\begin{equation*}\begin{split}
\mathcal{\hat{K}}\left[\kappa+\sum\limits_{i=1}^{n_1-1}\kappa_iv_i(x_1)+\sum\limits_{j=1}^{n_2-1}\delta_j\zeta_j(x_2)\right]= &\kappa+ \sum\limits_{i=1}^{n_1-1}\Psi_i(\kappa_1,\dots\kappa_{n_1-1},\delta_1,\dots,\delta_{n_2-1}) v_i(x_1)\\ +&\sum\limits_{j=1}^{n_2-1}\Psi_{n_1+j-1}(\kappa_1,\dots\kappa_{n_1-1},\delta_1,\dots,\delta_{n_2-1} )\zeta_j(x_2).
\end{split}
\end{equation*}
This gives the exact solution of the given equation \eqref{1.1} in the form
\[ u(x_1,x_2,t)= \Phi(t)+\sum\limits_{i=1}^{n_1-1}\Phi_i(t) v_i(x_1)+\sum\limits_{j=1}^{n_2-1}\Phi_{n_1+j-1}(t )\zeta_j(x_2).\]
Thus, the invariance conditions of type II linear spaces for the given differential operator $\mathcal{\hat{K}}[u]$ take the form
\begin{eqnarray*}
&&\mathcal{D}_{x_{1}}^{n_1}[\mathcal{\hat{K}}]\equiv \frac{d^{n_1}\mathcal{\hat{K}}}{dx_1^{n_1}}+a_{n_1-1}\frac{d^{n_1-1}\mathcal{\hat{K}}}{dx_1^{n_1-1}}\dots+a_1\frac{d\mathcal{\hat{K}}}{dx_1}=0,\\
&&\mathcal{D}_{x_2}^{n_2}[\mathcal{\hat{K}}]\equiv \frac{d^{n_2}\mathcal{\hat{K}}}{dx_2^{n_2}}+b_{n_2-1}\frac{d^{n_2-1}\mathcal{\hat{K}}}{dx_2^{n_2-1}}\dots+b_1\frac{d\mathcal{\hat{K}}}{dx_2}=0,\\
&&\frac{\partial^{2} \mathcal{\hat{K}}}{\partial x_1\partial x_2}=0
\end{eqnarray*}
along with
\begin{eqnarray*}
\begin{aligned}
&\mathcal{D}_{x_1}^{n_1}[u]\equiv \frac{d^{n_1}u}{dx_1^{n_1}}+a_{n_1-1}\frac{d^{n_1-1}u}{dx_1^{n_1-1}}+\dots+a_1\frac{d u}{dx_1}=0,\\
&\mathcal{D}_{x_2}^{n_2}[u]\equiv \frac{d^{n_2}u}{dx_2^{n_2}}+b_{n_2-1}\frac{d^{n_2-1}u}{dx_2^{n_2-1}}+\dots+b_1\frac{d u}{dx_2}=0,\\
&\frac{\partial^{2} u}{\partial x_1\partial x_2}=0
\end{aligned}
\end{eqnarray*}
and their differential consequences with respect to $x_1$ and $x_2$ if the differential operator $\mathcal{\hat{K}}[u] $ admits the type II linear space $\mathcal{V}_{n_1+n_2-1}$.

We would like to point out that the type II linear space $\mathcal{V}_{n_1+n_2-1}$ is the subspace of type I linear space $\mathcal{V}_{n_1n_2}$, that is, $\mathcal{V}_{n_1+n_2-1}\subseteq\mathcal{V}_{n_1n_2}$. Additionally, we observe that type II linear space $\mathcal{V}_{n_1+n_2-1}$ is smallest non-trivial (both variables involving $x_1$ and $x_2$) invariant subspace for the given differential operator $\mathcal{\hat{K}}[u]$.\\
For example, let us consider the type I and type II linear spaces $\mathcal{V}_{4}=Span\{1,x_1,x_2,x_1x_2\}$ and  $\mathcal{V}_{3}=Span\{1,x_1,x_2\}$  and also, the differential operator is \begin{eqnarray*}
&\mathcal{\hat{K}}[u]=& \dfrac{\partial }{\partial x_1}\left[(c_2u^2+c_1u+c_0)\left(\dfrac{\partial u}{\partial x_1}\right)\right]
+\dfrac{\partial }{\partial x_2}\left[(\beta_2u^2+\beta_1 u+\beta_0)\left(\dfrac{\partial u}{\partial x_2}\right)\right] \\&&
+({d_1 u+d_0}) \left(\dfrac{\partial u}{\partial x_1}\right)+(\lambda_1u+\lambda_0)\left(\dfrac{\partial u}{\partial x_2}\right)+k_1u+k_0,
\end{eqnarray*}
where $c_2,\beta_2, c_i,d_i,a_1,b_1,\beta_i,\lambda_i,k_i\in\mathbb{R},\ i=0,1$.  Note that $\mathcal{V}_{3}\subsetneqq\mathcal{V}_{4}$. Clearly, we check that the type I linear space is not invariant under $\mathcal{\hat{K}}[u]$ because
\begin{eqnarray*}\begin{aligned}
&\mathcal{\hat{K}}\left[\delta_1+\delta_2x_1+\delta_3x_2+\delta_4x_1x_2 \right]&\\
=&(k_1\delta_2+\delta_4\lambda_0+2c_2\delta_2^3+d_1\delta_2^2+
4\beta_2\delta_1\delta_3\delta_4+2\beta_2\delta_2\delta_3^2+
2\beta_1\delta_3\delta_4+\delta_1\delta_4\lambda_1+\delta_2\delta_3\lambda_1)
x_1&
\\&
+(k_1\delta_3+\delta_3^2\lambda_1+2\beta_2\delta_3^3+
4c_2\delta_1\delta_2\delta_4+2c_2\delta_2^2\delta_3+
2c_1\delta_2\delta_4+d_1\delta_1\delta_4+d_1\delta_2\delta_3+
d_0\delta_4)x_2
\\&
+(6\beta_2\delta_3^2\delta_4+6c_2\delta_2^2\delta_4+
2d_1\delta_2\delta_4+
2\delta_3\delta_4\lambda_1+k_1\delta_4  )x_1x_2+
(6c_2x_1\delta_2\delta_4^2+4c_2\delta_2\delta_3\delta_4\\&
+2c_2x_2\delta_3\delta_4^2+
2c_2\delta_1\delta_4^2+d_1\delta_3\delta_4+2c_2x_1x_2\delta_4^3+
d_1x_1\delta_4^2+\beta_1\delta_4^2)x_2^2+
(6x_2\beta_2\delta_3\delta_4^2\\&
+2x_1\beta_2\delta_2\delta_4^2+
2\beta_2\delta_1\delta_4^2+\delta_2\delta_4\lambda_1+2x_1x_2\beta_2\delta_4^3+
x_2\delta_4^2\lambda_1+c_1\delta_4^2+
4\beta_2\delta_2\delta_3\delta_4)x_1^2
\\
+&c_1\delta_2^2+\beta_1\delta_3^2+d_0\delta_2+\delta_3\lambda_0+k_1\delta_1+d_1\delta_1\delta_2
+
2\beta_2\delta_1\delta_3^2+\delta_1\delta_3\lambda_1
+2c_2\delta_1\delta_2^2+k_0,
\end{aligned}
\end{eqnarray*}
which gives $\mathcal{\hat{K}}\left[\delta_1+\delta_2x_1+\delta_3x_2+\delta_4x_1x_2 \right] \notin\mathcal{V}_{4}$.
But, the type II linear space $\mathcal{V}_{3}=Span\{1,x_1,x_2\}$ is invariant under the given differential operator $\mathcal{\hat{K}}[u]$, since
\begin{eqnarray*}\begin{aligned}
\mathcal{\hat{K}}\left[\delta_1+\delta_2x_1+\delta_3x_2\right]& = (2c_2\delta_2^3+d_1\delta_2^2+2\beta_2\delta_2\delta_3^2+
\delta_2\delta_3\lambda_1+k_1\delta_2)x_1+ (2c_2\delta_2^2\delta_3+\delta_3^2\lambda_1\\&
+2\beta_2\delta_3^3+d_1\delta_2\delta_3+k_1\delta_3)x_2+
2c_2\delta_1\delta_2^2+
2\beta_2\delta_1\delta_3^2+c_1\delta_2^2+
d_1\delta_1\delta_2\\
+&
\beta_1\delta_3^2+\delta_1\delta_3\lambda_1+d_0\delta_2+k_1\delta_1+
\delta_3\lambda_0+k_0 \in\mathcal{V}_{3}.
\end{aligned}
\end{eqnarray*}
From this, we obtain smallest non-trivial invariant subspace $\mathcal{V}_{3}$ which is proper subspace of $\mathcal{V}_{4}$.
Next, we would like to explain how to find the invariant subspaces for the given two-dimensional time-fractional generalized convection-reaction-diffusion-wave equation \eqref{1}.
\end{itemize}
\section{Invariant subspaces for the given equation \eqref{1} and its special types}
In this section, we explain how to construct the invariant subspaces for the generalized two-dimensional time-fractional convection-reaction-diffusion-wave equation \eqref{1} and its special kind of equations such as convection-diffusion-wave equation, reaction-diffusion-wave equation, and diffusion-wave equation.
\subsection{Invariant subspaces for the two-dimensional time-fractional\\ convection-reaction-diffusion-wave equation \eqref{1}}
In this subsection, we give a detailed systematic study for finding the invariant subspaces of the given equation \eqref{1}.
Thus, the differential operator $\mathcal{\hat{K}}[u]$ can be considered for the given equation \eqref{1} as follows
\begin{equation}\label{2.1}
\mathcal{\hat{K}}[u]=
 \sum\limits^{2}_{i=1} \frac{\partial }{\partial x_i}\left(A_i(u)\frac{\partial u}{\partial x_i}\right) + \sum\limits^{2}_{i=1}B_i(u)\frac{\partial u}{\partial x_i}+C(u).
\end{equation}
Here, we can discuss various dimensions of invariant subspaces for the above-mentioned differential operator $\mathcal{\hat{K}}[u]$.
 For example, in this work, we consider the following possibilities of $(n_1,n_2)$ as
 \begin{equation*}
 \left\{(1,1),(1,2),(1,3),(2,1),(2,2),(2,3),(3,1),(3,2),(3,3)\right\}.
 \end{equation*}
 From this, let us first consider $n_1=n_2=2$. Thus, we present a detailed study for finding the invariant subspaces of the given differential operator \eqref{2.1}.\\
\textbf{Estimation of invariant subspaces for the given equation \eqref{1}:} Here we would like to point out that two-types of linear spaces are available for the given differential operator \eqref{2.1}. Now, we define the linear spaces as
\begin{eqnarray*}
\begin{aligned}
&\mathcal{V}_{2}=\left\{y_1\ \big{|}\ \mathcal{D}^2_{x_1}[y_1]\equiv \frac{d^2y_1}{dx_1^{2}}+a_{1}\frac{dy_1}{dx_1}+a_{0}y_1=0\right\},\\
&\mathcal{V}_{2}=\left\{y_2\ \big{|}\ \mathcal{D}^2_{x_2}[y_2]\equiv \frac{d^{2}y_2}{dx_2^{2}}+b_{1}\frac{dy_2}{dx_2}+b_{0}y_2=0\right\},
\end{aligned}
\end{eqnarray*}
where the functions $y_i$, $i=1,2$, are the linear combinations of the two-linearly independent solutions of $\mathcal{D}^2_{x_i}[y_i]=0$. Then, let us assume that the functions $v_1(x_1)$ and $v_2(x_2)$ are two-linearly independent solutions of $\mathcal{D}^2_{x_1}[y_1]=0$. Similarly, we assume that another set of functions $\zeta_1(x_2)$ and $\zeta_2(x_2)$ are two-linearly independent solutions of $\mathcal{D}^2_{x_2}[y_2]=0$. Then we can define two-types of linear spaces for the given differential operator \eqref{2.1} that are discussed below. For this case, we obtain the dimensions of type I and type II invariant subspaces as follows
 \begin{itemize}
 \item[1.] $dim(\mathcal{V}_{n_1n_2})=4$ if either $1\notin\mathcal{V}_{n_1}$ and $1\notin\mathcal{V}_{n_2}$ or $1\in\mathcal{V}_{n_1} (or\ 1\in\mathcal{V}_{n_2})$,
 \item[2.] $dim(\mathcal{V}_{n_1+n_2-1})=3$ if $1\in\mathcal{V}_{n_1}$ $\&$ $1\in\mathcal{V}_{n_2}$,
 \end{itemize}
 which are discussed below.\\
\textbf{Type I and Type II linear spaces:}
First, we consider the type I linear space in the following form
\begin{eqnarray}
\begin{aligned}\label{2.2}
\mathcal{V}_{2\times2}&=\mathcal{V}_{4}=\text{Span}\left\{v_1(x_1)\zeta_1(x_2),v_1(x_1)\zeta_{2}(x_2),v_{2}(x_1)\zeta_1(x_2),v_{2}(x_1)\zeta_{2}(x_2)\right\}\\
&=\left\{\xi(x_1,x_2)=\sum\limits_{j=1}^{2}\left(\sum\limits_{i=1}^{2}c_{ij}v_i(x_1)\right)\zeta_j(x_2)\ \big{|}\ \mathcal{D}^2_{x_1}[\xi]=0, \mathcal{D}^2_{x_2}[\xi]=0\right\}
\end{aligned}
\end{eqnarray}
and type II linear space is considered as
\begin{eqnarray}
\begin{aligned}\label{2.3}
\mathcal{V}_{2+2-1}=&\mathcal{V}_{3}=\text{Span}\left\{1,v_1(x_1),\zeta_1(x_2)\right\}\\
=&\left\{\xi(x_1,x_2)=\kappa+\kappa_1v_1(x_1)+\delta_1\zeta_1(x_2)\ \big{|}\ \mathcal{D}^2_{x_k}[\xi]=0,\frac{\partial^2\xi}{\partial x_1\partial x_2}=0,\right.\\
&\left.a_0=b_0=0,k=1,2\right\}.
\end{aligned}
\end{eqnarray}
Now we wish to find the type I and type II invariance conditions for the differential operator $\mathcal{\hat{K}}[u]$ given in \eqref{2.1}.\\
\textbf{Type I invariant conditions:} If the given differential operator $\mathcal{\hat{K}}[u]$ admits type I linear space $\mathcal{V}_{4}$, then the type I invariant conditions of $\mathcal{\hat{K}}[u]$ read in the form
\begin{eqnarray*}
&&\mathcal{D}^2_{x_{1}}[\mathcal{\hat{K}}[u]]\equiv \frac{d^{2}\mathcal{\hat{K}}[u]}{dx_1^{2}}+a_{1}\frac{d\mathcal{\hat{K}}[u]}{dx_1}+a_{0}\mathcal{\hat{K}}[u]=0,\\
&&\mathcal{D}^2_{x_2}[\mathcal{\hat{K}}[u]]\equiv \frac{d^2\mathcal{\hat{K}}[u]}{dx_2^2}+b_1\frac{d\mathcal{\hat{K}}[u]}{dx_2}+b_{0}\mathcal{\hat{K}}[u]=0
\end{eqnarray*}
along with
\begin{eqnarray*}
\begin{aligned}
&\mathcal{D}^2_{x_1}[u]\equiv \frac{d^2u}{dx_1^2}+a_1\frac{d u}{dx_1}+a_{0}u=0,\\
&\mathcal{D}^2_{x_2}[u]\equiv \frac{d^2u}{dx_2^2}+b_1\frac{d u}{dx_2}+b_{0}u=0,
\end{aligned}
\end{eqnarray*}
where $ a_i,\ (i=0,1,2)$ and $b_j, \  (j=0,1,2)$ are constants to be determined, and the differential operator $\mathcal{\hat{K}}[u]$ given in \eqref{2.1}.
In a similar manner, we can find the type II invariant conditions for the given differential operator $\mathcal{\hat{K}}[u]$. \\
\textbf{Type II invariant conditions:} If the given differential operator $\mathcal{\hat{K}}[u]$ admits type II linear space $\mathcal{V}_{3}$, then the type II invariant conditions of $\mathcal{\hat{K}}[u]$ take the form
\begin{eqnarray}
\begin{aligned}\label{2.4}
&\mathcal{D}^2_{x_{1}}[\mathcal{\hat{K}}[u]]\equiv \frac{d^{2}\mathcal{\hat{K}}[u]}{dx_1^{2}}+a_{1}\frac{d\mathcal{\hat{K}}[u]}{dx_1}=0, \\ &\mathcal{D}^2_{x_2}[\mathcal{\hat{K}}[u]]\equiv \frac{d^{2}\mathcal{\hat{K}}[u]}{dx_2^{2}}+b_{1}\frac{d\mathcal{\hat{K}}[u]}{dx_2}=0,\\
&\frac{\partial^{2} \mathcal{\hat{K}}[u]}{\partial x_1\partial x_2}=0
 \end{aligned}
\end{eqnarray}
along with
\begin{eqnarray}
\begin{aligned}\label{2.5}
&\mathcal{D}^2_{x_{1}}[u]\equiv \frac{d^{2}u}{dx_1^2}+a_{1}\frac{d u}{dx_1}=0,\\
&\mathcal{D}^2_{x_2}[u]\equiv \frac{d^{2}u}{dx_2^{2}}+b_{1}\frac{du}{dx_2}=0,\\
&\frac{\partial^{2} u}{\partial x_1\partial x_2}=0,
 \end{aligned}
 \end{eqnarray}
 where $ a_i,\ (i=1,2)$ and $b_j, \  (j=1,2)$ are constants to be determined, and the differential operator $\mathcal{\hat{K}}[u]$ given in \eqref{2.1}.

Now, we give a detailed explanation for finding the type II invariant subspaces $\mathcal{V}_{3}$ of the given differential operator $\mathcal{\hat{K}}[u]$ using the above-discussed case $n_1=n_2=2$.
Thus, we substitute the given differential operator $\mathcal{\hat{K}}[u]$ into \eqref{2.4}, which gives
\begin{eqnarray}
\begin{aligned}\label{2.6}
&2(A_2)_{u}u_{x_{1}x_{2}}^2+4(A_2)_{uu}u_{x_{2}}u_{x_1}u_{x_{1}x_{2}}+[-2(A_2)_{u}b_1
+2(B_2)_{u}]u_{x_1}u_{x_{1}x_{2}}+(A_1)_{uuu}u_{x_{1}}^4\\
&+[-5(A_1)_{uu}a_1+(B_1)_{uu}]uA_{x_{1}}^3+u_{x_{2}}^2(A_2)_{uuu}u_{x_{1}}^2+[-(A_2)_{uu}b_1+(B_2)_{uu}]u_{x_{2}}u_{x_{1}}^2\\
&+[4(A_1)_{u}a_1^{2}-2(B_1)_{u}a_1+C_{uu}]u_{x_{1}}^2=0,
 \end{aligned}
 \end{eqnarray}
\begin{eqnarray}
\begin{aligned}\label{2.7}
&[-A_1a_1+B_1]u_{x_{1}x_{2}x_{2}}+2(A_1)_{u}u_{x_{1}}u_{x_{1}x_{2}x_{2}}+2(A_1)_{u}u_{x_{1}x_{2}}^2 +[(4(A_1)_{uu}u_{x_{2}}+
2(A_1)_{u}b_1)u_{x_{1}}\\
&+(-2(A_1)_{u}a_1+2(B_1)_{u})u_{x_{2}}-A_1b_1a_1+B_1b_1] u_{x_{1}x_{2}}+ (A_1)_{uuu}u_{x_{2}}^2u_{x_{1}}^2+(A_2)_{uuu}u_{x_{2}}^4\\&
+[-(A_1)_{uu}a_1+(B_1)_{uu}]u_{x_{2}}^2u_{x_{1}} + [-5(A_2)_{uu}b_1+(B_2)_{uu}]u_{x_{2}}^3\\
&+[4(A_2)_{u}b_1^{2}-2(B_2)_{u}b_1+C_{uu}]u_{x_{2}}^2
=0 \;
\end{aligned}
\end{eqnarray}
and
\begin{eqnarray}
 \begin{aligned}\label{2.8}
 &[A_2b_1+B_2+2(A_2)_{u}u_{x_2}]u_{x_2x_2x_1}+3(A_1)_{uu}u_{x_1}^{2}u_{x_1x_2}+[-6(A_1)_ua_1+2(B_1)_u]u_{x_1}u_{x_1x_2}\\
 &+3(A_2)_{uu}u_{x_2}^2u_{x_1x_2}+[2(B_2)_u-4(A_2)_ub_1]u_{x_2}u_{x_1x_2}+[A_1a_1^2-B_1a_1+C_u]u_{x_1x_2}\\
 &+(A_1)_{uuu}u_{x_{1}}^3u_{x_{2}}+[(B_1)_{uu}-3(A_1)_{uu}a_1]u_{x_{1}}^2u_{x_{2}}+
[(B_2)_{uu}-3(A_2)_{uu}a_1]u_{x_{2}}^2u_{x_{1}}\\
&+(A_2)_{uuu}u_{x_{2}}^3u_{x_{1}}+
[(A_1)_ua_{1}^{2}+(A_2)_ub_{1}^{2}-(B_1)_ua_{1}-(B_2)_ub_{1}+C_{uu}]u_{x_{1}}u_{x_{2}}=0.
\end{aligned}
\end{eqnarray}
Substituting equation \eqref{2.5} into equations \eqref{2.6}-\eqref{2.8}, which are reduce to the following form
 \begin{eqnarray}
 \begin{aligned} \label{2.9}
& (A_1)_{uuu}u_{x_{1}}^4+
[-5(A_1)_{uu}a_1+(B_1)_{uu}]u_{x_{1}}^3+u_{x_{2}}^2(A_2)_{uuu}u_{x_{1}}^2\\
&+[-(A_2)_{uu}b_1+(B_2)_{uu}]u_{x_{2}}u_{x_{1}}^2+[4(A_1)_{u}a_1^{2}-2(B_1)_{u}a_1+C_{uu}]u_{x_{1}}^2=0,
 \end{aligned}
 \end{eqnarray}
\begin{eqnarray}
 \begin{aligned}\label{2.10}
& (A_1)_{uuu}u_{x_{2}}^2u_{x_{1}}^2+(A_2)_{uuu}u_{x_{2}}^4+
[-(A_1)_{uu}a_1+(B_1)_{uu}]u_{x_{2}}^2u_{x_{1}}  \\
&+[-5(A_2)_{uu}b_1+(B_2)_{uu}]u_{x_{2}}^3+[4(A_2)_{u}b_1^{2}-2(B_2)_{u}b_1+C_{uu}]u_{x_{2}}^2=0
\end{aligned}
 \end{eqnarray}
 {and}
\begin{eqnarray}
 \begin{aligned}\label{2.11}
&(A_1)_{uuu}u_{x_{1}}^3u_{x_{2}}+[(B_1)_{uu}-3(A_1)_{uu}a_1]u_{x_{1}}^2u_{x_{2}}+
[(B_2)_{uu}-3(A_2)_{uu}a_1]u_{x_{2}}^2u_{x_{1}}\\
&+(A_2)_{uuu}u_{x_{2}}^3u_{x_{1}}+
[(A_1)_ua_{1}^{2}+(A_2)_ub_{1}^{2}-(B_1)_ua_{1}-(B_2)_ub_{1}+C_{uu}]u_{x_{1}}u_{x_{2}}=0,
\end{aligned}
 \end{eqnarray}
 where $(A_i)_u=\dfrac{dA_i}{du}$, $(A_i)_{uu}=\dfrac{d^2A_i}{du^2}$,$(B_i)_u=\dfrac{dB_i}{du}$, $(B_i)_{uu}=\dfrac{d^2B_i}{du^2}$, $u_{x_i}=\dfrac{\partial u}{\partial x_i}$, $i=1,2$ and etc.
The above equations \eqref{2.9}-\eqref{2.11} are reduces to the system of over-determined equations. In general, the obtained system of over-determined equations may not be solvable. Thus, let us first consider the functions $A_1(u)=c_2u^2+c_1u+c_0$, $A_2(u)=\beta_2 u^2+\beta_1 u+\beta_0$, $B_1(u)=d_2u^2+d_1u+d_0$, $B_2(u)=\lambda_2 u^2+\lambda_1u+\lambda_0$ and $C(u)= k_3u^3+k_2u^2+k_1u+k_0$. \\
For this case, we obtain the cubic non-linear differential operator $\mathcal{\hat{K}}[u]$ that is obtained as
\begin{equation}\begin{split}
\mathcal{\hat{K}}[u]\equiv & \frac{\partial }{\partial x_1}\left[(c_2u^2+c_1u+c_0)\left(\frac{\partial u}{\partial x_1}\right)\right] +
\frac{\partial }{\partial x_2}\left[(\beta_2 u^2+\beta_1 u+\beta_0)\left(\frac{\partial u}{\partial x_2}\right)\right]\\ & +(d_2u^2+d_1u+d_0) \left(\frac{\partial u}{\partial x_1}\right)+(\lambda_2 u^2+\lambda_1u+\lambda_0)\left(\frac{\partial u}{\partial x_2}\right)\\ &+ k_3u^3+k_2u^2+k_1u+k_0,\label{2.12}\end{split}
\end{equation}
where $c_i,\beta_i,d_i,\lambda_i,k_i,k_3\in\mathbb{R}$, $i=0,1,2$.
Substituting the functions $A_i(u)$, $B_i(u)$ $(i=1,2)$ and $C(u)$ into the above equations \eqref{2.9}-\eqref{2.11}, which yields the following system of over-determined equations
\begin{eqnarray*}
    &6a_1c_2-2d_2 = 0 ,
   6b_1\beta_2-2\lambda_2 = 0,
     2d_2 -2a_1c_2 = 0,
  2\lambda_2-10b_1\beta_2 = 0,
   2d_2 -10a_1c_2 = 0,\\&
8a_1^2c_2-4a_1d_2+6k_3 = 0,
 4a_1^2c_1-2a_1d_1+2k_2 = 0,
 8b_1^2\beta_2-4b_1\lambda_2+6k_3 = 0 ,\\&
 2\lambda_2 -2b_1\beta_2 = 0, 4b_1^2\beta_1-2b_1\lambda_1+2k_2 = 0,
  a_1^2c_1+b_1^2\beta_1-a_1d_1-b_1\lambda_1+2k_2 = 0,\\&
  a_1^2c_2+b_1^2\beta_2-a_1d_2-b_1\lambda_2+3k_3 = 0.
 \end{eqnarray*}
Solving the above over-determined system of equations, we get type II three-dimensional invariant subspaces $\mathcal{V}_{3}$ with their corresponding differential operators that are given in Table 2 and Table 3. Proceeding the above similar procedure, we can find various dimensions of the invariant subspaces for the given differential operator \eqref{2.1} using the other possible values of $(n_1,n_2)$ with different non-linearities that are listed from Table 4 to Table 15.
\begin{sidewaystable}
\caption{Discussion for various dimensions of the invariant subspaces for the differential operator $\mathcal{\hat{K}}[u]$ given in \eqref{2.1}}
\begin{center}
\begin{tabular}{ccccc}
\toprule
Cases & Differential operator \eqref{2.1}  &
 Different types of  & Reference\\
& with different non-linearities& invariant subspaces  &     tables \\
\midrule

1.& Quadratic  non-linearities
  &
           $Type II:\mathcal{V}_3=\mathcal{V}_{2+2-1}$

 &
          Tables 2,4-7 \\

2.& Cubic  non-linearities
  & $\left\{
         \begin{array}{ll}
           Type I:\mathcal{V}_1=\mathcal{V}_{1\cdot1}\\
           Type II:\mathcal{V}_3=\mathcal{V}_{2+2-1}\\
         \end{array}
       \right.
$
 &  Tables 3, 7-12\\
3.& Other  non-linearities
  &
          $ Type I:\mathcal{V}_2=\mathcal{V}_{1\cdot2}$

 &
          Tables 13-15
           \\

\bottomrule
\end{tabular}
\end{center}
\label{0}

\caption{Quadratic differential operator \eqref{2.12} and their corresponding invariant subspace}
\begin{center}
\begin{tabular}{cccc}
\toprule
Cases & $\mathcal{\hat{K}}[u]=\sum\limits^{2}_{i=1}(A_i(u)u_{x_{i}})_{x_{i}}+\sum\limits^{2}_{i=1}B_i(u)u_{x_{i}}+C(u)$  &  Invariant\\
 &  $A_i(u),\ B_i(u)\ \&\ C(u),i=1,2$ &     Subspaces \\
\midrule
1.& $\left\{
         \begin{array}{ll}
         \left\{
         \begin{array}{ll}
            A_1(u)=\frac{-b_1^{2}\beta_1}{a_1^{2}}u+c_{0} \\
            A_2(u)=\beta_1u+\beta_{0}\\
            B_1(u)=\frac{-4b_1^{2}\beta_1+b_1\lambda_1}{a_1}u+d_{0}\\
            B_2(u)=\lambda_{1}u+\lambda_{0}\\
            C(u)=({-2b_{1}^{2}\beta_1+b_1\lambda_1})u^{2}+k_{1}u+k_{0}
         \end{array}
       \right.\\
       \left\{
         \begin{array}{ll}
A_1(u)=c_{0} \\
            A_2(u)=\frac{-c_1a_1^{2}}{b_1^2}u+\beta_{0}\\
            B_1(u)=\frac{4a_1^{2}c_1+b_1\lambda_1}{a_1}u+d_{0}\\
            B_2(u)=\lambda_{1}u+\lambda_{0}\\
            C(u)=({2a_{1}^{2}c_1+b_1\lambda_1})u^{2}+k_{1}u+k_{0}\\ \end{array}
       \right. \end{array}
       \right.

$  & $\left.
                                                                                                                                   \begin{array}{ll}

                                                                                                                                     \mathcal{V}_{3}=Span\{1,e^{-a_1x_1},e^{-b_1x_2}\}&
                                                                                                                                   \end{array}
                                                                                                                                 \right.
$\\
\bottomrule
\end{tabular}
\end{center}
\label{default}
\end{sidewaystable}
\begin{sidewaystable}
\caption{Cubic differential operator \eqref{2.12} and their corresponding invariant subspaces}
\begin{center}
\begin{tabular}{cccc}
\toprule
Cases & $\mathcal{\hat{K}}[u]=\sum\limits^{2}_{i=1}(A_i(u)u_{x_{i}})_{x_{i}}+\sum\limits^{2}_{i=1}B_i(u)u_{x_{i}}+C(u)$  &  Invariant\\
 &  $A_i(u),\ B_i(u)\ \&\ C(u),i=1,2$ &     Subspaces \\
\midrule
1.&$ \left\{
         \begin{array}{ll}
           A_1(u)=c_{2}u^{2}+c_{1}u+c_{0}\\
          A_2(u)=\beta_{2}u^{2}+\beta_{1}u+\beta_{0}\quad \qquad \qquad\\
          B_1(u)=d_{1}u+d_{0}\\
          B_2(u)=\lambda_{1}u+\lambda_{0}\\
          C(u)=k_{1}u+k_{0}

         \end{array}
       \right.

$  & $\left.
                                                                                                                                   \begin{array}{ll}

                                                                                                                                     \mathcal{V}_3=Span\{1,x_1,x_2\} &
                                                                                                                                   \end{array}
                                                                                                                                 \right.
$\\
\bottomrule
\end{tabular}
\end{center}
\label{default}
\end{sidewaystable}

\begin{sidewaystable}
\caption{Differential operator \eqref{2.1} with quadratic non-linearities and their corresponding invariant subspaces}
\begin{center}
\begin{tabular}{cccc}
\toprule
Cases & $\mathcal{\hat{K}}[u]=\sum\limits^{2}_{i=1}(A_i(u)u_{x_{i}})_{x_{i}}+\sum\limits^{2}_{i=1}B_i(u)u_{x_{i}}+C(u)$  &  Invariant\\
 &  $A_i(u),\ B_i(u)\ \&\ C(u),i=1,2$ &     Subspaces \\
\midrule

1. & $\left\{
         \begin{array}{ll}
         \left\{
         \begin{array}{ll}
           A_1(u)=c_{1}u+c_{0}\\
          A_2(u)=\dfrac{-c_{1}a_{1}^{2}}{b_{1}^{2}}u+\beta_{0}\\
          B_1(u)=d_{1}u+d_{0}\\
          B_2(u)=\frac{-4a_{1}^{2}c_{1}+a_{1}d_{1}}{b_{1}}u+\lambda_{0}\\
          C(u)={(-2a_{1}^{2}c_{1}+a_{1}d_{1})}u^{2}+k_{1}u+k_{0}
         \end{array}
       \right.\\
       \left\{
         \begin{array}{ll}
            A_1(u)=c_{0} \\
            A_2(u)=\beta_{0}\\
            B_1(u)=\frac{b_1\lambda_1}{a_1}u+d_{0}\\
            B_2(u)=\lambda_{1}u+\lambda_{0}\\
            C(u)=b_1\lambda_1u^{2}+k_{1}u+k_{0}\\
         \end{array}
       \right.
       \end{array}
       \right.
$  & $\left.
                                                                                                                                   \begin{array}{ll}

                                                                                                                                     \mathcal{V}_{3}=Span\{1,e^{-a_1x_1},e^{-b_1x_2}\}&
                                                                                                                                   \end{array}
                                                                                                                                 \right.
$\\
\bottomrule
\end{tabular}
\end{center}
\label{default}
\end{sidewaystable}

\begin{sidewaystable}
\caption{Differential operator \eqref{2.1} with quadratic non-linearities and their corresponding invariant subspaces}
\begin{center}
\begin{tabular}{cccc}
\toprule
Cases & $\mathcal{\hat{K}}[u]=\sum\limits^{2}_{i=1}(A_i(u)u_{x_{i}})_{x_{i}}+\sum\limits^{2}_{i=1}B_i(u)u_{x_{i}}+C(u)$  &  Invariant\\
 &  $A_i(u),\ B_i(u)\ \&\ C(u),i=1,2$ &     Subspaces \\
\midrule
4. & $\left\{
         \begin{array}{ll}
           A_1(u)=c_{1}u+c_{0}\\
          A_2(u)=\beta_{1}u+\beta_{0}\\
          B_1(u)=d_{0}\\
          B_2(u)=\lambda_{1}u+\lambda_{0} \\
          C(u)=k_{1}u+k_{0}
         \end{array}
       \right.
$  & $\left.
                                                                                                                                   \begin{array}{ll}

                                                                                                                                     \mathcal{V}_{4}=Span\{1,x_1,\frac{x_1^{2}}{2},x_2\}&
                                                                                                                                   \end{array}
                                                                                                                                 \right.
$\\
\bottomrule
\end{tabular}
\end{center}
\label{default}
\end{sidewaystable}

\begin{sidewaystable}
\caption{Differential operator \eqref{2.1} with quadratic non-linearities and their corresponding invariant subspaces}
\begin{center}
\begin{tabular}{cccc}
\toprule
Cases & $\mathcal{\hat{K}}[u]=\sum\limits^{2}_{i=1}(A_i(u)u_{x_{i}})_{x_{i}}+\sum\limits^{2}_{i=1}B_i(u)u_{x_{i}}+C(u)$  &  Invariant\\
 &  $A_i(u),\ B_i(u)\ \&\ C(u),i=1,2$ &     Subspaces \\
\midrule
8. & $\left\{
         \begin{array}{ll}
           A_1(u)=c_{1}u+c_{0}\qquad\qquad\\
          A_2(u)=\beta_{0}\\
          B_1(u)=d_{0}\\
          B_2(u)=\lambda_{0}\\
          C(u)=k_{1}u+k_{0}
         \end{array}
       \right.
$  & $\left.
                                                                                                                                   \begin{array}{ll}

                                                                                                                                     \mathcal{V}_5=Span\{1,e^{\frac{1}{2}(-b_2+\sqrt{b_{2}^{2}-4b_{1}})x_2}, e^{-\frac{1}{2}(b_2+\sqrt{b_{2}^{2}-4b_{1}})x_2},x_1,\frac{{x_1}^{2}}{2}\}&
                                                                                                                                   \end{array}
                                                                                                                                 \right.
$\\
\bottomrule
\end{tabular}
\end{center}
\label{default}
\end{sidewaystable}

\begin{sidewaystable}
\caption{Differential operator \eqref{2.1} with quadratic non-linearities and their corresponding invariant subspaces}
\begin{center}
\begin{tabular}{cccc}
\toprule
Cases & $\mathcal{\hat{K}}[u]=\sum\limits^{2}_{i=1}(A_i(u)u_{x_{i}})_{x_{i}}+\sum\limits^{2}_{i=1}B_i(u)u_{x_{i}}+C(u)$  &  Invariant\\
 &  $A_i(u),\ B_i(u)\ \&\ C(u),i=1,2$ &     Subspaces \\
\midrule
12. & $\left\{
         \begin{array}{ll}
           A_1(u)=c_{1}u+c_{0}\\
          A_2(u)=\beta_{0}\\
          B_1(u)=d_{0}\\
          B_2(u)=\lambda_{0}\\
          C(u)=k_{1}u+k_{0}
         \end{array}
       \right.
$  & $\left.
                                                                                                                                   \begin{array}{ll}

                                                                                                                                     \mathcal{V}_{5}=Span\{1,x_1,\frac{x_1^{2}}{2},sin{(\sqrt{b_1}x_2)}, cos{(\sqrt{b_1}x_2)}\}&
                                                                                                                                   \end{array}
                                                                                                                                 \right.
$\\
\bottomrule
\end{tabular}
\end{center}
\label{default}

\caption{Differential operator \eqref{2.1} with cubic non-linearities and their corresponding invariant subspaces}
\begin{center}
\begin{tabular}{cccc}
\toprule
Cases & $\mathcal{\hat{K}}[u]=\sum\limits^{2}_{i=1}(A_i(u)u_{x_{i}})_{x_{i}}+\sum\limits^{2}_{i=1}B_i(u)u_{x_{i}}+C(u)$  &  Invariant\\
 &  $A_i(u),\ B_i(u)\ \&\ C(u),i=1,2$ &     Subspaces \\
\midrule

 1. & $\left\{
         \begin{array}{ll}
           A_1(u)=( \frac{-3b_{0}^{2}\beta_{2}+a_{0}d_{2}+b_{0}\lambda_{2}}{3a_{0}^{2}})u^{2}+ ( \frac{-2b_{0}^{2}\beta_{1}+a_{0}d_{1}+b_{0}\lambda_{1}-k_{2}}{
       2a_{0}^{2}})u+c_{0}\\
          A_2(u)=\beta_{2}u^{2}+\beta_{1}u+\beta_{0}\\
          B_1(u)=d_{2}u^{2}+d_{1}u+d_{0}\\
          B_2(u)=\lambda_{2}u^{2}+\lambda_{1}u+\lambda_{0}\\
          C(u)=k_{2}u^{2}+k_{1}u
         \end{array}
       \right.
$  & $\left.
                                                                                                                                   \begin{array}{ll}

                                                                                                                                     \mathcal{V}_1=Span\{e^{-a_0x_1-b_0x_2}\}&
                                                                                                                                   \end{array}
                                                                                                                                 \right.
$\\
\bottomrule
\end{tabular}
\end{center}
\label{default}
\end{sidewaystable}
\begin{sidewaystable}
\caption{Differential operator \eqref{2.1} with cubic non-linearities and their corresponding invariant subspaces}
\begin{center}
\begin{tabular}{cccc}
\toprule
Cases & $\mathcal{\hat{K}}[u]=\sum\limits^{2}_{i=1}(A_i(u)u_{x_{i}})_{x_{i}}+\sum\limits^{2}_{i=1}B_i(u)u_{x_{i}}+C(u)$  &  Invariant\\
 &   $A_i(u),\ B_i(u)\ \&\ C(u),i=1,2 $ &     Subspaces \\
\midrule
2. & $\left\{
         \begin{array}{ll}

           A_1(u)=c_{0}\\
          A_2(u)=\dfrac{\lambda_2}{3b_0}u^2+\beta_1u+\beta_{0}\\
          B_1(u)=d_{0}\\
          B_2(u)=\lambda_2u^2+\lambda_1u+\lambda_{0}\\
          C(u)=(-2b_{0}^2\beta_{1}+b_0\lambda_1)u^{2}+k_{1}u
         \end{array}
       \right.
$  & $\left.
                                                                                                                                   \begin{array}{ll}
                                                                                                                                   \mathcal{V}_2 = Span\{e^{-b_{0}x_2},e^{-a_{1}x_1-b_0x_2}\}&
\end{array}
       \right.$\\
\bottomrule
\end{tabular}
\end{center}
\label{default}
\end{sidewaystable}
\begin{sidewaystable}
\caption{Differential operator \eqref{2.1} with cubic non-linearities and their corresponding invariant subspaces}
\begin{center}
\begin{tabular}{cccc}
\toprule
Cases & $\mathcal{\hat{K}}[u]=\sum\limits^{2}_{i=1}(A_i(u)u_{x_{i}})_{x_{i}}+\sum\limits^{2}_{i=1}B_i(u)u_{x_{i}}+C(u)$  &  Invaraint\\
 &   $A_i(u),\ B_i(u)\ \&\ C(u),\ i=1,2$ &     Subspaces \\
\midrule
5. & $\left\{
         \begin{array}{ll}
           A_1(u)=c_{2}u^2+c_{0}\\
          A_2(u)=\beta_{0}\\
          B_1(u)=3a_{0}c_{2}u^2+d_{1}u+d_{0}\qquad\\
          B_2(u)=\lambda_{0}\\
        C(u)=a_{0}d_{1}u^2+k_{1}u\\
         \end{array}
       \right.
$  & $\left.
                                                                                                                                   \begin{array}{ll}

                                                                                                                                     \mathcal{V}_{2}=Span\{e^{-a_{0}x_1}sin{(\sqrt{b_{0}}x_2)},e^{-a_{0}x_1}cos{(\sqrt{b_{0}}x_2)}\}&
                                                                                                                                   \end{array}
                                                                                                                                 \right.
$\\
\bottomrule
\end{tabular}
\end{center}
\label{default}
\end{sidewaystable}
\begin{sidewaystable}
\caption{Differential operator \eqref{2.1} with cubic non-linearities and their corresponding invariant subspaces}
\begin{center}
\begin{tabular}{cccc}
\toprule
Cases & $\mathcal{\hat{K}}[u]=\sum\limits^{2}_{i=1}(A_i(u)u_{x_{i}})_{x_{i}}+\sum\limits^{2}_{i=1}B_i(u)u_{x_{i}}+C(u)$  &  Invariant\\
 &  $A_i(u),\ B_i(u)\ \&\ C(u),i=1,2$ &     Subspaces \\
\midrule

9. & $\left\{
         \begin{array}{ll}
           A_1(u)=c_{1}u+c_{0}\\
          A_2(u)=\beta_{2}u^{2}+\beta_{1}u+\beta_{0}\\
          B_1(u)=d_{0}\\
          B_2(u)=\lambda_{1}u+\lambda_{0}\\
          C(u)=k_{1}u+k_{0}
         \end{array}
       \right.
$  & $\left.
                                                                                                                                   \begin{array}{ll}

                                                                                                                                     \mathcal{V}_4=Span\{1,x_1,\frac{x_{1}^{2}}{2},x_2\} &
                                                                                                                                   \end{array}
                                                                                                                                 \right.
$\\
\bottomrule
\end{tabular}
\end{center}
\label{default}
\end{sidewaystable}

\begin{sidewaystable}
\caption{Differential operator \eqref{2.1} with cubic non-linearities and their corresponding invariant subspaces}
\begin{center}
\begin{tabular}{cccc}
\toprule
Cases & $\mathcal{\hat{K}}[u]=\sum\limits^{2}_{i=1}(A_i(u)u_{x_{i}})_{x_{i}}+\sum\limits^{2}_{i=1}B_i(u)u_{x_{i}}+C(u)$  &  Invariant\\
 &  $A_i(u),\ B_i(u)\ \&\ C(u),i=1,2$ &     Subspaces \\
\midrule

13. & $\left\{
         \begin{array}{ll}
           A_1(u)=c_{2}u^2+c_{1}u+c_{0}\\
          A_2(u)=\beta_{0}\\
          B_1(u)=d_{1}u+d_{0}\\
          B_2(u)=\lambda_{0} \\
          C(u)=k_{1}u+k_{0}
         \end{array}
       \right.
$  & $\left.
                                                                                                                                   \begin{array}{ll}

                                                                                                                                     \mathcal{V}_4=Span\{1,x_{1},sin(\sqrt{b_{1}}x_{2}),cos(\sqrt{b_{1}}x_{2})\}&

                                                                                                                                   \end{array}
                                                                                                                                 \right.
$\\
\bottomrule
\end{tabular}
\end{center}
\label{default}

\caption{Differential operator \eqref{2.1} with other non-linearities and their corresponding invariant subspaces}
\begin{center}
\begin{tabular}{cccc}
\toprule
Cases & $\mathcal{\hat{K}}[u]=\sum\limits^{2}_{i=1}(A_i(u)u_{x_{i}})_{x_{i}}+\sum\limits^{2}_{i=1}B_i(u)u_{x_{i}}+C(u)$  &  Invariant\\
 &  $A_i(u), B_i(u)\ \&\ C(u),i=1,2$ &     Subspaces \\
 \midrule

1. & $\left\{
         \begin{array}{ll}
          A_1(u)=\frac{d_{3}}{4a_{0}}u^{3}+c_{2}u^{2}+c_{0}\\
          A_2(u)=\beta_{0}\\
          B_1(u)=d_{3}u^{3}+d_{2}u^{2}+d_{1}u+d_{0}\\
          B_2(u)=\lambda_{0}\\
          C(u)=(-3a_{0}^{2}c_{2}+a_{0}d_{2})u^{3}+ (a_{0}d_{1})u^{2}+k_{1}u
         \end{array}
       \right.
$  & $\left.
                                                                                                                                   \begin{array}{ll}

                                                                                                                                     \mathcal{V}_2=Span\{
                                                                                                                                     e^{-a_{0}x_1-\frac{1}{2}(b_1-\sqrt{b_{1}^{2}-4b_{0}})x_2},
                                                                                                                                  e^{-a_{0}x_1-\frac{1}{2}(b_1+\sqrt{b_{1}^{2}-4b_{0}})x_2}
                                                                                                                                  \}&
                                                                                                                                   \end{array}
                                                                                                                                 \right.
$\\
\bottomrule
\end{tabular}
\end{center}
\label{default}
\end{sidewaystable}
\begin{sidewaystable}
\caption{Differential operator \eqref{2.1} with other non-linearities and their corresponding invariant subspaces}
\begin{center}
\begin{tabular}{cccc}
\toprule
Cases & $\mathcal{\hat{K}}[u]=\sum\limits^{2}_{i=1}(A_i(u)u_{x_{i}})_{x_{i}}+\sum\limits^{2}_{i=1}B_i(u)u_{x_{i}}+C(u)$  &  Invariant\\
 &  $A_i(u), B_i(u)\ \&\ C(u),i=1,2$ &     Subspaces \\
 \midrule
2. & $\left\{
         \begin{array}{ll}
           A_1(u)=c_{0}\\
          A_2(u)=\frac{1}{4}\frac{\lambda_{3}}{b_{0}}u^{3}+\beta_{1}u+\beta_{0}\\
          B_1(u)=d_{0}\\
          B_2(u)=\lambda_{3}u^{3}+\lambda_{2}u^{2}+\lambda_{1}u+\lambda_{0}\\
          C(u)=b_{0}\lambda_{2}u^{3}+(-2b_{0}^{2}\beta_{1}+ b_{0}\lambda_{1})u^{2}\\\qquad+k_{1}u\\
         \end{array}
       \right.
$  & $\left.
                                                                                                                                   \begin{array}{ll}

                                                                                                                                     \mathcal{V}_2=Span\{e^{-b_{0}x_2-\frac{1}{2}(a_1-\sqrt{a_{1}^{2}-4a_{0}})x_1},  e^{-b_{0}x_2-\frac{1}{2}(a_1+\sqrt{a_{1}^{2}-4a_{0}})x_1}\}&
                                                                                                                                 \end{array}
                                                                                                                                 \right.
$\\
\bottomrule
\end{tabular}
\end{center}
\label{default}
\end{sidewaystable}

\begin{sidewaystable}
\caption{Differential operator \eqref{2.1} with other non-linearities and their corresponding invariant subspaces}
\begin{center}
\begin{tabular}{cccc}
\toprule
Cases & $\mathcal{\hat{K}}[u]=\sum\limits^{2}_{i=1}(A_i(u)u_{x_{i}})_{x_{i}}+\sum\limits^{2}_{i=1}B_i(u)u_{x_{i}}+C(u)$  &  Invariant\\
 &  $A_i(u), B_i(u)\ \&\ C(u),i=1,2$ &     Subspaces \\
 \midrule

5. & $\left\{
         \begin{array}{ll}
         A_1(u)=c_0\\
          A_2(u)=\beta_3 u^3+\beta_1 u+\beta_0\\
          B_1(u)=d_0\\
          B_2(u)=4b_0 \beta_3 u^3+\lambda_2 u^2+\lambda_1 u+\lambda_0\\
          C(u)=\lambda_2 b_0 u^3 +(-2 b_0^2 \beta_1 +b_0 \lambda_1)u^2+k_1u\\

         \end{array}
       \right.
$  & $\left.
                                                                                                                                   \begin{array}{ll}

                                                                                                                                  \mathcal{V}_{2}=Span\{e^{-b_0x_2}sin({\sqrt{a_{0}}x_1}),e^{-b_0x_2}cos{(\sqrt{a_{0}}x_1)}\}&
                                                                                                                                   \end{array}
                                                                                                                                 \right.
$\\
\bottomrule
\end{tabular}
\end{center}
\label{default}
\end{sidewaystable}

\newpage
\subsection{Special kinds of the given equation \eqref{1} with their invariant subspaces}
In this subsection, we show how to find the invariant subspaces for the special types of the given two-dimensional convection-reaction-diffusion-wave equation \eqref{1}. It is interesting to note that the given equation \eqref{1} is reduced into three different special kinds of equations that are discussed below.
\begin{itemize}
\item[$\bullet$] When $C(u)=0$, the given equation \eqref{1} becomes the generalised two-dimensional time-fractional convection-diffusion-wave equation
\begin{equation}\label{4.1}
\frac{\partial^{\alpha}u}{\partial t^{\alpha}}=
 \sum\limits^{2}_{i=1} \frac{\partial }{\partial x_i}\left(A_i(u)\frac{\partial u}{\partial x_i}\right) + \sum\limits^{2}_{i=1}B_i(u)\frac{\partial u}{\partial x_i},\ \alpha\in(0,2].
\end{equation}
\item[$\bullet$] If $B_1(u)=B_2(u)=0$, then the given equation \eqref{1} can be viewed as the generalised two-dimensional time-fractional reaction-diffusion-wave equation
\begin{equation} \label{4.2}
\frac{\partial^{\alpha}u}{\partial t^{\alpha}}=
 \sum\limits^{2}_{i=1} \frac{\partial }{\partial x_i}\left(A_i(u)\frac{\partial u}{\partial x_i}\right)+C(u),\ \alpha\in(0,2].
\end{equation}
\item[$\bullet$] If $B_1(u)=B_2(u)=C(u)=0$, then the equation \eqref{1} reduces into
\begin{equation}
\label{4.3}
\frac{\partial^{\alpha}u}{\partial t^{\alpha}}=
 \sum\limits^{2}_{i=1} \frac{\partial }{\partial x_i}\left(A_i(u)\frac{\partial u}{\partial x_i}\right), \alpha\in(0,2]
\end{equation}
which is familiarly known as the generalized two-dimensional time-fractional diffusion-wave equation.
\end{itemize}
Proceeding the above similar procedure, we can find the various dimensions of the type I and type II invariant subspaces for the above-mentioned equations with different power-law non-linearities. The obtained invariant subspaces with their corresponding differential operators are listed from Table 17 to Table 27.
\begin{table}
\caption{Discussion of Type I and Type II invariant subspaces for the  equations \eqref{4.1}-\eqref{4.3} with different power-law non-linearities}
\begin{center}
\begin{tabular}{ccccc}
\toprule
Cases& Differential operators  &
 Type of invariant & Reference\\
& for the equations \eqref{4.1}-\eqref{4.3} &  Subspaces  &     Tables \\
\midrule
1.& $\mathcal{\hat{K}}[u]=\sum\limits^{2}_{i=1}(A_i(u)u_{x_{i}})_{x_{i}}+\sum\limits^{2}_{i=1}B_i(u)u_{x_{i}}$
  & $\left\{
         \begin{array}{ll}
           Type I:\mathcal{V}_1=\mathcal{V}_{1\cdot1}\\
           Type II:\mathcal{V}_3=\mathcal{V}_{2+2-1}\\
         \end{array}
       \right.
$
 & Tables 17 to 23 \\
2.&$\mathcal{\hat{K}}[u]=\sum\limits^{2}_{i=1}(A_i(u)u_{x_{i}})_{x_{i}}+C(u)$ & $\left\{
         \begin{array}{ll}
           Type I:\mathcal{V}_1=\mathcal{V}_{1\cdot1}\\
           Type II:\mathcal{V}_4=\mathcal{V}_{2+3-1}\\
         \end{array}
       \right.
$  & Tables 24, 25\\
3.&  $\mathcal{\hat{K}}[u]=\sum\limits^{2}_{i=1}(A_i(u)u_{x_{i}})_{x_{i}}$
  & $\left\{
         \begin{array}{ll}
           Type I:\mathcal{V}_1=\mathcal{V}_{1\cdot1}\\
           Type II:\mathcal{V}_3=\mathcal{V}_{2+2-1}\\
         \end{array}
       \right.
$  & Tables 26,27 \\
\bottomrule
\end{tabular}
\end{center}
\label{default}
\end{table}
\begin{table}
\caption{Differential operator of \eqref{4.1} with quadratic non-linearities and their corresponding invariant subspaces}
\begin{center}
\begin{tabular}{cccc}
\toprule
Cases &
 $\mathcal{\hat{K}}[u]=\sum\limits^{2}_{i=1}(A_i(u)u_{x_{i}})_{x_{i}}+\sum\limits^{2}_{i=1}B_i(u)u_{x_{i}}$  &  Invariant\\
 &  $A_i(u)\ \&\ B_i(u),i=1,2 $ &     Subspaces \\
\midrule

1.& $\left\{
         \begin{array}{ll}
           A_1(u)=\dfrac{-b_{1}\lambda_{1}}{2a_{1}^{2}}u+c_{0}\\
          A_2(u)=\dfrac{\lambda_{1}}{2b_{1}}u+\beta_{0}\\
          B_1(u)=\frac{-b_{1}\lambda_{1}}{a_{1}}u+d_{0}\\
          B_2(u)=\lambda_{1}u+\lambda_{0}\\

         \end{array}
       \right.
$  & $\left.
                                                                                                                                   \begin{array}{ll}

                                                                                                                                     \mathcal{V}_{3}=Span\{1,e^{-a_1x_1},e^{-b_1x_2}\}&
                                                                                                                                   \end{array}
                                                                                                                                 \right.
$\\
\bottomrule
\end{tabular}
\end{center}
\label{default}
\end{table}
                                                                                                                                   \begin{table}
\caption{Differential operator of \eqref{4.1} with quadratic non-linearities and their corresponding invariant subspaces}
\begin{center}
\begin{tabular}{cccc}
\toprule
Cases &
 $\mathcal{\hat{K}}[u]=\sum\limits^{2}_{i=1}(A_i(u)u_{x_{i}})_{x_{i}}+\sum\limits^{2}_{i=1}B_i(u)u_{x_{i}}$  &  Invariant\\
 &   $A_i(u)\ \&\ B_i(u),i=1,2$ &     Subspaces \\
\bottomrule
\end{tabular}
\end{center}
\label{default}
\end{table}
\begin{table}
\caption{Differential operator of \eqref{4.1} with cubic non-linearities and their corresponding invariant subspaces}
\begin{center}
\begin{tabular}{cccc}
\toprule
Cases &  $\mathcal{\hat{K}}[u]=\sum\limits^{2}_{i=1}(A_i(u)u_{x_{i}})_{x_{i}}+\sum\limits^{2}_{i=1}B_i(u)u_{x_{i}}$  &  Invariant\\
 &   $A_i(u)\ \&\ B_i(u),i=1,2$ &     Subspaces \\
\midrule

1. & $\left\{
         \begin{array}{ll}
           A_1(u)=\frac{1}{3} \frac{-3b_{0}^{2}\beta_{2}+a_{0}d_{2}+b_{0}\lambda_{2}}{a_{0}^{2}}u^{2}
           \\\quad    \qquad +\frac{1}{2}\frac{-2b_{0}^{2}\beta_{1}+a_{0}d_{1}+b_{0}\lambda_{1}}{a_{0}^{2}}u+c_{0}\\
          A_2(u)=\beta_{2}u^{2}+\beta_{1}u+\beta_{0}\\
          B_1(u)=d_{2}u^{2}+d_{1}u+d_{0}\\
          B_2(u)=\lambda_{2}u^{2}+\lambda_{1}u+\lambda_{0}\\

         \end{array}
       \right.
$  & $\left.
                                                                                                                                   \begin{array}{ll}

                                                                                                                                     \mathcal{V}_{1}=Span\{e^{-a_0x_1-b_0x_2}\}&
                                                                                                                                   \end{array}
                                                                                                                                 \right.
$\\
\bottomrule
\end{tabular}
\end{center}
\label{default}
\end{table}
\begin{table}
\caption{Differential operator of \eqref{4.1} with cubic non-linearities and their corresponding invariant subspaces}
\begin{center}
\begin{tabular}{cccc}
\toprule
Cases &  $\mathcal{\hat{K}}[u]=\sum\limits^{2}_{i=1}(A_i(u)u_{x_{i}})_{x_{i}}+\sum\limits^{2}_{i=1}B_i(u)u_{x_{i}}$ &  Invariant\\
 &  $A_i(u)\ \&\ B_i(u),i=1,2 $ &     Subspaces \\
\midrule

6. & $\left\{
         \begin{array}{ll}
           A_1(u)=c_{2}u^{2}+c_{1}u+c_{0}\\
          A_2(u)=\beta_{0}\\
          B_1(u)=d_{1}u+d_{0}\\
          B_2(u)=\lambda_{0}\\

         \end{array}
       \right.
$  & $\left.
                                                                                                                                   \begin{array}{ll}

                                                                                                                                     \mathcal{V}_{3}=Span\{1,e^{-b_1x_2},x_1\}&
                                                                                                                                   \end{array}
                                                                                                                                 \right.
$\\
\bottomrule
\end{tabular}
\end{center}
\label{default}
\end{table}
\begin{table}
\caption{Differential operator of \eqref{4.1} with cubic non-linearities and their corresponding invariant subspaces}
\begin{center}
\begin{tabular}{cccc}
\toprule
Cases &  $\mathcal{\hat{K}}[u]=\sum\limits^{2}_{i=1}(A_i(u)u_{x_{i}})_{x_{i}}+\sum\limits^{2}_{i=1}B_i(u)u_{x_{i}}$  &  Invariant\\
 &   $A_i(u)\ \&\ B_i(u), i=1,2$    & Subspaces \\
\midrule

12. & $\left\{
         \begin{array}{ll}
           A_1(u)=c_{0}\\
          A_2(u)=\beta_2 u^2 +\beta_1 u+\beta_{0}\\
          B_1(u)=d_{0}\\
          B_2(u)=\lambda_1 u+\lambda_{0}\\

         \end{array}
       \right.
$  & $\left.
                                                                                                                                   \begin{array}{ll}

                                                                                                                                     \mathcal{V}_{4}=Span\{1,x_2,sin(\sqrt{a_{1}}x_1),\\
                                                                                                                                     cos(\sqrt{a_{1}}x_1)\}&
                                                                                                                                   \end{array}
                                                                                                                                 \right.
$\\
\bottomrule
\end{tabular}
\end{center}
\label{default}
\end{table}

\begin{table}
\caption{Differential operator of \eqref{4.1} with other non-linearities and their corresponding invariant subspaces}
\begin{center}
\begin{tabular}{cccc}
\toprule
Cases & $\mathcal{\hat{K}}[u]=\sum\limits^{2}_{i=1}(A_i(u)u_{x_{i}})_{x_{i}}+\sum\limits^{2}_{i=1}B_i(u)u_{x_{i}}$ &  Invariant\\
 &  $A_i(u)\ \&\ B_i(u),i=1,2 $ &     Subspaces \\
\midrule

1. & $\left\{
         \begin{array}{ll}
           A_1(u)=\dfrac{d_{3}}{4a_{0}}u^{3}+\dfrac{d_{2}}{3a_{0}}u^{2}+c_{0}\\
          A_2(u)=\beta_{0}\\
          B_1(u)=d_{3}u^{3}+d_{2}u^{2}+d_{0}\\
          B_2(u)=\lambda_{0}\\

         \end{array}
       \right.
$  & $\left.
                                                                                                                                   \begin{array}{ll}

                                                                                                                                     \mathcal{V}_{2}=Span\{e^{-a_{0}x_1-\frac{1}{2}(b_1-\sqrt{b_{1}^{2}-4b_{0}})x_2},\\ e^{-a_{0}x_1-\frac{1}{2}(b_1+\sqrt{b_{1}^{2}-4b_{0}})x_2}\}&
                                                                                                                                   \end{array}
                                                                                                                                 \right.
$\\
\bottomrule
\end{tabular}
\end{center}
\label{default}
\end{table}
\begin{table}
\caption{Differential operator of \eqref{4.1} with other non-linearities and their corresponding invariant subspaces}
\begin{center}
\begin{tabular}{cccc}
\toprule
Cases & $\mathcal{\hat{K}}[u]=\sum\limits^{2}_{i=1}(A_i(u)u_{x_{i}})_{x_{i}}+\sum\limits^{2}_{i=1}B_i(u)u_{x_{i}}$  &  Invariant\\
 &  $A_i(u)\ \&\ B_i(u),i=1,2$ &     Subspaces \\
\midrule
3. & $\left\{
         \begin{array}{ll}
           A_1(u)=\dfrac{d_{3}}{4a_{0}}u^{3}+\dfrac{d_{2}}{3a_{0}}u^{2}\\\qquad\quad
           +\frac{d_{1}}{2a_{0}}u + c_{0}\\
          A_2(u)=\beta_{0}\\
          B_1(u)=d_{3}u^{3}+d_{2}u^{2}+d_{1}u+d_{0}\\
          B_2(u)=\lambda_{0}\\

         \end{array}
       \right.
$  & $\left.
                                                                                                                                   \begin{array}{ll}

                                                                                                                                     \mathcal{V}_{3}=Span\{e^{-a_{0}x_1-\frac{1}{2}(b_1-\sqrt{b_{1}^{2}-4b_{0}})x_2},\\ e^{-a_{0}x_1-\frac{1}{2}(b_1+\sqrt{b_{1}^{2}-4b_{0}})x_2},e^{-a_{0}x_1}\}&
                                                                                                                                   \end{array}
                                                                                                                                 \right.
$\\
\bottomrule
\end{tabular}
\end{center}
\label{default}
\end{table}

\begin{table}
\caption{Differential operator of \eqref{4.2} with quadratic non-linearities and their corresponding invariant subspaces}
\begin{center}
\begin{tabular}{cccc}
\toprule
Cases &  $\mathcal{\hat{K}}[u]=\sum\limits^{2}_{i=1}(A_i(u)u_{x_{i}})_{x_{i}}+C(u)$  &  Invariant\\
 &  $A_i(u)\ \&\ C(u),i=1,2 $ &     Subspaces \\
\midrule

1. & $\left\{
         \begin{array}{ll}
           A_1(u)=\frac{-1}{2} \frac{2b_{0}^{2}\beta_{1}+k_{2}}{a_{0}^{2}}u+c_{0}\\
          A_2(u)=\beta_{1}u+\beta_{0}\\
          C(u)=k_{2}u^{2}+k_{1}u\\

         \end{array}
       \right.
$  & $\left.
                                                                                                                                   \begin{array}{ll}

                                                                                                                                     \mathcal{V}_{1}=Span\{e^{-a_0x_1-b_0x_2}\}&
                                                                                                                                   \end{array}
                                                                                                                                 \right.
$\\
\bottomrule
\end{tabular}
\end{center}
\label{default}
\end{table}
\begin{table}
\caption{Differential operator of \eqref{4.2} with cubic non-linearities and their corresponding invariant subspaces}
\begin{center}
\begin{tabular}{cccc}
\toprule
Cases &  $\mathcal{\hat{K}}[u]=\sum\limits^{2}_{i=1}(A_i(u)u_{x_{i}})_{x_{i}}+C(u)$ &  Invariant\\
 &  $A_i(u)\ \&\ C(u), i=1,2 $ &     Subspaces \\
\midrule
1. & $\left\{
         \begin{array}{ll}
           A_1(u)=c_2 u^2+c_1 u+c_{0}\\
          A_2(u)=\beta_{0}\\
          C(u)=k_{1}u+k_0\\

         \end{array}
       \right.
$  & $\left.
                                                                                                                                   \begin{array}{ll}

                                                                                                                                     \mathcal{V}_{3}=Span\{1,x_1,e^{-b_{1}x_{2}}\}&
                                                                                                                                   \end{array}
                                                                                                                                 \right.
$\\
2. & $\left\{
         \begin{array}{ll}
           A_1(u)=c_{0}\\
          A_2(u)=\beta_2 u^2+\beta_1 u+\beta_{0}\\
          C(u)=k_{1}u+k_0\\

         \end{array}
       \right.
$  & $\left.
                                                                                                                                   \begin{array}{ll}

                                                                                                                                     \mathcal{V}_{3}=Span\{1,x_2,e^{-a_{1}x_{1}}\}&
                                                                                                                                   \end{array}
                                                                                                                                 \right.
$\\
\bottomrule
\end{tabular}
\end{center}
\label{default}
\end{table}

\begin{sidewaystable}
\caption{Differential operator of \eqref{4.3} with quadratic non-linearities and their corresponding invariant subspaces}
\begin{center}
\begin{tabular}{cccc}
\toprule
Cases & $\mathcal{\hat{K}}[u]=\sum\limits^{2}_{i=1}(A_i(u)u_{x_{i}})_{x_{i}}$ &  Invariant\\
 &  $A_i(u),i=1,2 $ &     Subspaces \\
\midrule
1.& $\left\{
         \begin{array}{ll}
              A_1(u)=c_{0}\\
          A_2(u)=\beta_{1}u+\beta_{0}\\

         \end{array}
       \right.
$  & $
\left\lbrace
                                                                                                                                   \begin{array}{ll}
                                                                                                                                                                                                                                                                    \mathcal{V}_{4}=Span\{1,e^{-a_{1}x_1},x_2,\frac{x_2^{2}}{2}\}\\

                                                                                                                                     \mathcal{V}_{4}=Span\{1,e^{-a_{2}x_1},x_1,x_2\}\\
                                                                                                                                    \mathcal{V}_{5}=Span\{1,e^{-\frac{1}{2}(a_{2}-\sqrt{a_{2}^{2}-4a_{1}})x_1},e^{-\frac{1}{2}(a_{2}+\sqrt{a_{2}^{2}-4a_{1}})x_1},x_2,\frac{x_2^{2}}{2}\}&

                                                                                                                                   \end{array}
                                                                                                                                                                                                                                                           \right.
$\\
\bottomrule
\end{tabular}
\end{center}
\label{default}
\end{sidewaystable}

\begin{sidewaystable}
\caption{Differential operator of \eqref{4.3} with cubic non-linearities and their corresponding invariant subspaces}
\begin{center}
\begin{tabular}{cccc}
\toprule
Cases & $\mathcal{\hat{K}}[u]=\sum\limits^{2}_{i=1}(A_i(u)u_{x_{i}})_{x_{i}}$ &  Invariant\\
 &  $A_i(u),i=1,2 $ &     Subspaces \\
\midrule
1.& $\left\{
         \begin{array}{ll}
           A_1(u)= \dfrac{-b_{0}^{2}\beta_{2}}{a_{0}^{2}}u^{2}+c_{0}\\
          A_2(u)=\beta_{2}u^{2}+\beta_{0}\\

         \end{array}
       \right.
$  & $\left.
                                                                                                                                   \begin{array}{ll}

                                                                                                                                     \mathcal{V}_{1}=Span\{e^{-a_0x_1-b_0x_2}\}&
                                                                                                                                   \end{array}
                                                                                                                                 \right.
$\\
\bottomrule
\end{tabular}
\end{center}
\label{default}
\end{sidewaystable}
\newpage
\section{Exact solutions of \eqref{1} along with initial conditions}
In this section, we would like to explain how to derive exact solutions for the initial value problem of the given two-dimensional time-fractional non-linear convection-reaction-diffusion-wave equation \eqref{1} with different types of non-linearities using the obtained invariant subspaces that are discussed in the previous section. Hence, let us first explain how to construct exact solutions for the initial value problem of the given equation \eqref{1} with cubic non-linearity using the obtained exponential subspaces.
\subsection{Exact solutions of \eqref{1} with cubic non-linearity using the exponential subspaces}
Now, we first discuss how to derive the exact solutions for the initial value problem of the two-dimensional time-fractional convection-reaction-diffusion-wave equation \eqref{1} with cubic non-linearity using the two-dimensional exponential subspace which was discussed in table 9 of case 3.
Thus, we consider the convection-reaction-diffusion-wave equation with cubic non-linearity in the form
\begin{eqnarray}
\begin{aligned}\label{5.1}
\frac{\partial^{\alpha}u}{\partial t^{\alpha}}=&
 \frac{\partial }{\partial x_1}\left[\left(\frac{d_2}{3a_0}u^2+c_1u+c_0\right)\frac{\partial u }{\partial x_1}\right]+\beta_0\frac{\partial^2 u }{\partial x_{2}^{2}}+(d_2u^2+d_1u+d_0)\frac{\partial u }{\partial x_1}\\
&+\lambda_0\frac{\partial u }{\partial x_2}+(-2a_0^2c_1+a_0d_1)u^2+k_1u,\ \alpha\in(0,2],\ t\geq0,\ x_1,x_2\in\mathbb{R}
\end{aligned}
\end{eqnarray}
along with the initial conditions
\begin{itemize}
\item[$\bullet$] If $\alpha\in(0,1]$, then $u(x_1,x_2,0)=\nu_1e^{-a_0x_1}+\nu_2e^{-(a_0x_1+b_1x_2)}$.
\item[$\bullet$] If $\alpha\in(1,2]$, then $u(x_1,x_2,0)=\nu_1e^{-a_0x_1}+\nu_2e^{-(a_0x_1+b_1x_2)}$ $\&$\\
${\dfrac{\partial u}{\partial t}}\mid_{t=0}=\mu_1e^{-a_0x_1}+\mu_2e^{-(a_0x_1+b_1x_2)}$.
\end{itemize}
For this case, we obtain the cubic non-linear differential operator
\begin{eqnarray*}
\mathcal{\hat{K}}[u]=&\left(\dfrac{d_2}{3a_0}u^2+c_1u+c_0\right)\dfrac{\partial^2 u}{\partial x_1^2}+\left(\dfrac{2d_2}{3a_0}u+c_1\right)\left(\dfrac{\partial u}{\partial x_1}\right)^2+\beta_0\dfrac{\partial^2 u }{\partial x_{2}^{2}}\\
&+(d_2u^2+d_1u+d_0)\dfrac{\partial u }{\partial x_1}
+\lambda_0\dfrac{\partial u }{\partial x_2}+(-2a_0^2c_1+a_0d_1)u^2+k_1u,
\end{eqnarray*}
which admits a two-dimensional exponential subspace $\mathcal{V}_2=Span\{e^{-a_0x_1},e^{-(a_0x_1+b_1x_2)}\}$, since for some constants $ \delta_1, \delta_2\in \mathbb{R} $
\begin{equation*}
\begin{split}
\mathcal{\hat{K}}[\delta_1 e^{-a_0x_1}+\delta_2e^{-(a_0x_1+b_1x_2)}]=&(a_0^2c_0-d_0a_0+k_1)\delta_1 e^{-a_0x_1}\\
&+(k_1+\beta_0b_1^2-\lambda_0b_1+c_0a_0^2-a_{0}d_{0})\delta_2e^{-(a_0x_1+b_1x_2)}\in\mathcal{V}_2.
\end{split}
\end{equation*}
Hence the equation \eqref{5.1} possesses an exact solution in the form
\begin{equation}\label{5.2}
u(x_1,x_2,t)=\Phi_1(t)e^{-a_0x_1}+\Phi_2(t)e^{-(a_0x_1+b_1x_2)},
\end{equation}
where the unknown functions $\Phi_1(t)$ and $\Phi_2(t)$ are to be determined.\\
Substituting \eqref{5.2} into \eqref{5.1}, we get the system of fractional ODEs
\begin{eqnarray*}
&&\label{5.17}\dfrac{d^{\alpha}\Phi_1(t)}{d t^{\alpha}}=(a_0^2c_0-d_0a_0+k_1){\Phi_1(t)}, \\
&&\label{5.18}\dfrac{d^{\alpha}\Phi_2(t)}{d t^{\alpha}}=(k_1+\beta_0b_1^2-\lambda_0b_1+c_0a_0^2-a_{0}d_{0}){\Phi_2(t)},\ \alpha\in(0,2].
\end{eqnarray*}
The above equations can be written into the following form
\begin{equation}\label{5.3}
\dfrac{d^{\alpha}\Phi_i(t)}{d t^{\alpha}}=\gamma_i{\Phi_i(t)}, i=1,2, \alpha\in(0,2],
\end{equation}
where $\gamma_1=(a_0^2c_0-d_0a_0+k_1)$ and $\gamma_2=(k_1+\beta_0b_1^2-\lambda_0b_1+c_0a_0^2-a_{0}d_{0})$.
First, let us consider $\alpha\in (0,1]$. It should be noted that the Laplace transformation of the Caputo fractional derivative of order $\alpha>0$ is given by \cite{kai,pi}
 \[ \mathcal{L}\left[\frac{d^\alpha w(\xi)}{d\xi^\alpha}\right]=s^\alpha\mathcal{L}[w(\xi)]-\sum_{m=0}^{n-1}s^{\alpha-(m+1)}w^{(m)}(0),\ n-1<\alpha\leq n,\ n\in\mathbb{N}, \quad Re(s)>0, \]
 where $w^{(m)}(0)= \dfrac{d^mw(\xi)}{d\xi}\mid_{\xi=0} ,\xi \in [0,\infty)$ and $\mathcal{L}[w(\xi)]=\tilde{w}(s)=\int\limits_{0}^{\infty}e^{-st}w(t)dt.$ Thus applying the Laplace transformation to equation \eqref{5.3}, we obtain
 \begin{equation*}
 s^\alpha\tilde{\Phi}_i(s)-s^{\alpha-1}\Phi_i(0)=\gamma_i\tilde{\Phi}_i(s).
 \end{equation*}
 Applying the inverse Laplace transformation for the above equation, we get
 \begin{equation*}
 \Phi_i(t)=\nu_i{E}_{\alpha,1}(\gamma_it^\alpha),\ \alpha\in(0,1],\ i=1,2,
 \end{equation*}
where $\nu_i=\Phi_i(0)$ and $E_{\alpha,\beta}(t^\alpha)$ is the two-parameter Mittag-Leffler function \cite{mathai}, defined as
$$E_{\alpha,\beta}(t^\alpha)=\sum\limits_{m=0}^{\infty}\dfrac{(t^\alpha)^m}{\Gamma(m\alpha+\beta)}.$$
Hence, for the case $\alpha\in(0,1]$, the obtained exact solution for the two-dimensional time-fractional cubic non-linear convection-reaction-diffusion equation \eqref{5.1} associated with two-dimensional exponential subspace $\mathcal{V}_2=Span\{e^{-a_0x_1},e^{-(a_0x_1+b_1x_2)}\}$ as
\begin{equation}\label{s1}
u(x_1,x_2,t)=\nu_1e^{-a_0x_1}{E}_{\alpha,1}(\gamma_1t^\alpha)+\nu_2e^{-(a_0x_1+b_1x_2)}{E}_{\alpha,1}(\gamma_2t^\alpha),\ \alpha\in(0,1],
\end{equation}
 where $\gamma_1=(a_0^2c_0-d_0a_0+k_1)$, $\gamma_2=(k_1+\beta_0b_1^2-\lambda_0b_1+c_0a_0^2-a_{0}d_{0})$ and $a_0,b_1,\nu_i\in\mathbb{R}$, $i=1,2$. Also, it should be noted that the obtained exact solution \eqref{s1} satisfies the initial condition $u(x_1,x_2,0)=\nu_1 e^{-a_0x_1}+\nu_2 e^{-(a_0x_1+b_1x_2)}$, since ${E}_{\alpha,1}(\gamma_it^\alpha)\mid_{ t=0}=1$.

Now, we consider the case $\alpha\in (1,2]$. Then, applying the Laplace transformation to equation \eqref{5.3}, we obtain
 \begin{equation*}
 \mathcal{L}\left[\dfrac{d^{\alpha}\Phi_{i}(t)}{d t^{\alpha}}\right]=\mathcal{L}\left[\gamma_{i}{\Phi_{i}(t)}\right], \alpha\in(1,2],
 \end{equation*}
 which gives
\begin{equation*}
s^\alpha\tilde{\Phi}_i(s)-s^{\alpha-1}\Phi_i(0)-s^{\alpha-2}\Phi_i'(0)=\gamma_i\tilde{\Phi}_i(s).
\end{equation*}
The above equation can be written as
\begin{equation*}
\tilde{\Phi}_i(s)=\Phi_i(0)\left(\dfrac{s^{\alpha-1}}{s^\alpha-\gamma_i}\right)+\Phi_i'(0)\left(\dfrac{s^{\alpha-2}}{s^\alpha-\gamma_i}\right).
\end{equation*}
Applying the inverse Laplace transformation of the above equation, we get
\begin{equation*}
 \Phi_i(t)=\nu_i E_{\alpha,1}(\gamma_it^{\alpha})+t\mu_iE_{\alpha,2}(\gamma_it^{\alpha}),\ \alpha\in(1,2],
\end{equation*}
where $\nu_i=\Phi_i(0)$ and $\mu_i=\Phi_{i}'(0)=\dfrac{d \Phi_i(t)}{dt}\big{|}_{ t=0}$, $i=1,2$.\\
For this case, we obtain an exact solution for the two-dimensional time-fractional convection-reaction-diffusion-wave equation with cubic non-linearity \eqref{5.1} as
\begin{eqnarray}
\begin{aligned}\label{5.5}
u(x_1,x_2,t)=&  (\nu_1E_{\alpha,1}(\gamma_1t^{\alpha})+t\mu_1E_{\alpha,2}(\gamma_1t^{\alpha}))e^{-a_0x_1}\\
&+(\nu_2E_{\alpha,1}(\gamma_2t^{\alpha})+t\mu_2E_{\alpha,2}(\gamma_2t^{\alpha}))e^{-(a_0x_1+b_1x_2)},\ \alpha\in (1,2],
\end{aligned}
\end{eqnarray}
 where $\gamma_1=(a_0^2c_0-d_0a_0+k_1)$, $\gamma_2=(k_1+\beta_0b_1^2-\lambda_0b_1+c_0a_0^2-a_{0}d_{0})$ and $a_0,b_1,\mu_i,\nu_i\in\mathbb{R}$, $i=1,2$.
Additionally, we observe that the obtained exact solution \eqref{5.5} satisfies the given initial conditions for $\alpha\in(1,2]$ as
\begin{equation*}
u(x_1,x_2,0)=\nu_1e^{-a_0x_1}+\nu_2e^{-(a_0x_1+b_1x_2)}
\end{equation*}
and
\begin{eqnarray*}
\begin{aligned}
\dfrac{\partial u}{\partial t}\mid_{t=0}=&\left(\nu_1t^{\alpha-1}E_{\alpha,\alpha}(\gamma_1t^\alpha)+\mu_1{E}_{\alpha,1}(\gamma_1t^\alpha)\right)e^{-a_0x_1}\\
&+\left(\nu_2t^{\alpha-1}E_{\alpha,\alpha}(\gamma_2t^\alpha)+\mu_2{E}_{\alpha,1}(\gamma_2t^\alpha)\right)e^{-(a_0x_1+b_1x_2)} \mid_{t=0}
\end{aligned}
\end{eqnarray*}
which gives
\begin{equation*}
\dfrac{\partial u}{\partial t}\mid_{t=0}=\mu_1e^{-a_0x_1}+\mu_2e^{-(a_0x_1+b_1x_2)}.
\end{equation*}

Next, we consider the initial value problem of two-dimensional time-fractional convection-reaction-diffusion-wave equation with cubic non-linearity
\begin{eqnarray}
\begin{aligned}\label{5.6}
\dfrac{\partial^{\alpha}u}{\partial t^{\alpha}}=&\dfrac{\partial }{\partial x_1}\left[\left(  \dfrac{-3b_{0}^{2}\beta_{2}+a_{0}d_{2}+b_{0}\lambda_{2}}{3a_{0}^{2}}u^{2}+ \dfrac{-2b_{0}^{2}\beta_{1}+a_{0}d_{1}+b_{0}\lambda_{1}-k_{2}}{2a_{0}^{2}}u+c_{0}\right)\dfrac{\partial u }{\partial x_1}\right]\\
&+\dfrac{\partial }{\partial x_2}\left[\left(\beta_{2}u^{2}+\beta_{1}u+\beta_{0}\right)\dfrac{\partial u }{\partial x_{2}}\right]+
(d_2u^2+d_1u+d_0)\dfrac{\partial u }{\partial x_1}\\
&+(\lambda_2u^2+\lambda_1u+\lambda_0)\dfrac{\partial u }{\partial x_2}+k_2u^2+k_1u,\ \alpha\in(0,2],\ t\geq0
\end{aligned}
\end{eqnarray}
subject to the initial conditions
\begin{eqnarray}
&&\label{i1} u(x_1,x_2,0)= \nu e^{-(a_0x_1+b_0x_2)} \  \text{if}\ \alpha \in(0,1].\\
&&\label{i2}u(x_1,x_2,0)= \nu e^{-(a_0x_1+b_0x_2)}\ \&
{\frac{\partial u}{\partial t}}\mid_{t=0}=\mu e^{-(a_0x_1+b_0x_2)} \ \textsl{if}\ \alpha \in(1,2].
\end{eqnarray}
Note that equation \eqref{5.6} admits a one-dimensional invariant subspace $\mathcal{V}_1=Span\{e^{-a_0x_1-b_0x_2}\}$ which was discussed in table 8 of case 1. Following the above similar procedure,
 the obtained exact solutions for the two-dimensional time-fractional convection-reaction-diffusion-wave equation with cubic non-linearity \eqref{5.6} are as follows
 \begin{eqnarray*}\label{5.7}
&& u(x_1,x_2,t)=\nu E_{\alpha,1}(\gamma t^\alpha)e^{-(a_0x_1+b_0x_2)}\ \text{if}\ \alpha \in (0,1],\\
&& u(x_1,x_2,t)=\left[\nu E_{\alpha,1}(\gamma t^\alpha)+t\mu E_{\alpha,2}(\gamma t^\alpha)\right]e^{-(a_0x_1+b_0x_2)}\ \text{if}\ \alpha\in (1,2],
 \end{eqnarray*}
where $\gamma=a_0^2c_0+b_0^2\beta_0-\lambda_0a_0-d_0b_0+k_1$ and $a_0,b_0,\nu,\mu\in\mathbb{R}$. In addition, note that the obtained exact solutions satisfy the given initial conditions \eqref{i1} and \eqref{i2}.

Finally, let us consider the two-dimensional time-fractional convection-reaction-diffusion-wave equation with cubic non-linearity
\begin{eqnarray}
\begin{aligned}\label{5.8}
\frac{\partial^{\alpha}u}{\partial t^{\alpha}}=&c_0\frac{\partial^2 u }{\partial x_{1}^{2}}+\frac{\partial }{\partial x_2}\left[(\frac{\lambda_2}{3b_0}u^2+\beta_1u+\beta_0)\frac{\partial u }{\partial x_2}\right]+d_0\frac{\partial u }{\partial x_1}\\
&+(\lambda_2u^2+\lambda_1u+\lambda_0)\frac{\partial u }{\partial x_2}+b_0\lambda_1u^2+k_1u,\ \alpha\in(0,2],\ t\geq0,
\end{aligned}
\end{eqnarray}
subject to the initial conditions
\begin{eqnarray}
&&\label{i3} \text{If}\ \alpha \in(0,1], \text{then}\; u(x_1,x_2,0)=\nu_1e^{-b_0x_2}+\nu_2e^{-(a_1x_1+b_0x_2)}.\\
&&\label{i4} \text{If}\ \alpha \in(1,2], \text{then} \left\{
              \begin{array}{ll}
                u(x_1,x_2,0)=\nu_1e^{-b_0x_2}+\nu_2e^{-(a_1x_1+b_0x_2)}\ \& & \\
                {\dfrac{\partial u}{\partial t}}\mid_{t=0} =\mu_1e^{-b_0x_2}+\mu_2e^{-(a_1x_1+b_0x_2)}. &
              \end{array}
            \right.
\end{eqnarray}
The above equation \eqref{5.8} admits an invariant subspace $\mathcal{V}_2=Span\{e^{-b_0x_2},e^{-(a_1x_1+b_0x_2)}\}$ that was discussed in case 2 of table 9.
In a similar manner, we can find exact solutions for the two-dimensional time-fractional convection-reaction-diffusion-wave equation \eqref{5.8} with cubic non-linearity in the form
 \begin{eqnarray*}
 \begin{aligned}\label{5.9}
u(x_1,x_2,t)=&\nu_1e^{-b_0x_2}{E}_{\alpha,1}(\gamma_1t^\alpha)+\nu_2e^{-(a_1x_1+b_0x_2)}{E}_{\alpha,1}(\gamma_2t^\alpha)\ \text{if}\ \alpha\in (0,1],\\
u(x_1,x_2,t)=&\left[\nu_1{E}_{\alpha,1}(\gamma_1t^\alpha)+t\mu_1{E}_{\alpha,2}(\gamma_1t^\alpha)\right]e^{-b_0x_2}\\
&+\left[\nu_2{E}_{\alpha,1}(\gamma_2t^\alpha)+t\mu_2{E}_{\alpha,2}(\gamma_2t^\alpha)\right]e^{-(a_1x_1+b_0x_2)}\ \text{if}\ \alpha\in (1,2],
\end{aligned}
 \end{eqnarray*}
where $\gamma_1=(b_0^2\beta_0-\lambda_0b_0+k_1),\gamma_2=(k_1+\beta_0b_0^2-\lambda_0b_0+c_0a_1^2-a_{1}d_{0})$ and $a_1,b_0,\nu_i, \mu_i\in\mathbb{R}$, $i=1,2$. Note that the above obtained exact solutions satisfy the given initial conditions \eqref{i3} and \eqref{i4}.

\subsection{Exact solutions of \eqref{1} with cubic non-linearity using the combination of exponential and trigonometric subspaces}
In this subsection, we consider the two-dimensional time-fractional convection-reaction-diffusion-wave equation \eqref{1} with cubic non-linearity in the form
\begin{eqnarray}
\begin{aligned}\label{5.10}
\frac{\partial^{\alpha}u}{\partial t^{\alpha}}=&\frac{\partial }{\partial x_1}\left[(c_2u^2+c_0)\frac{\partial u }{\partial x_1}\right]+\beta_0\frac{\partial^2 u }{\partial x_{2}^{2}}+(3a_0c_2u^2+d_1u+d_0)\frac{\partial u }{\partial x_1}\\
&+a_0d_1u^2+k_1u,\ \alpha\in(0,2],\ t\geq0,
\end{aligned}
\end{eqnarray}
subject to the initial conditions
\begin{eqnarray}
&&\label{ic1} \text{If}\ \alpha \in(0,1], \text{then}\; u(x_1,x_2,0)= e^{-a_0x_1}\left(\nu_1sin(\sqrt{b_0}x_2)+\nu_2cos(\sqrt{b_0}x_2)\right).\\
&&\label{ic2}\text{If}\ \alpha \in(1,2], \text{then} \left\{
              \begin{array}{ll}
                u(x_1,x_2,0)=e^{-a_0x_1}\left(\nu_1sin(\sqrt{b_0}x_2)+\nu_2cos(\sqrt{b_0}x_2)\right)\& \\
                {\dfrac{\partial u}{\partial t}}\mid_{t=0} = e^{-a_0x_1}\left(\mu_1sin(\sqrt{b_0}x_2)+\mu_2cos(\sqrt{b_0}x_2)\right). &
              \end{array}
            \right.
\end{eqnarray}
Equation \eqref{5.10} admits an invariant subspace $\mathcal{V}_2=Span\{e^{-a_0x_1}sin(\sqrt{b_0}x_2),e^{-a_0x_1}cos(\sqrt{b_0}x_2)\}$ that is discussed in case 5 of table 10. For this case, we can find exact solutions for the the two-dimensional time-fractional convection-reaction-diffusion-wave equation with cubic non-linearity \eqref{5.10} in the form
\begin{eqnarray}
 \begin{aligned}\label{5.11}
 u(x_1,x_2,t)=&e^{-a_0x_1}E_{\alpha,1}(\gamma t^\alpha)\left[\nu_1sin(\sqrt{b_0}x_2)+\nu_2cos(\sqrt{b_0}x_2)\right]\ \text{if}\ \alpha \in (0,1],  \\
 u(x_1,x_2,t)=&e^{-a_0x_1}E_{\alpha,1}(\gamma t^\alpha)\left[\nu_1sin(\sqrt{b_0}x_2)+\nu_1cos(\sqrt{b_0}x_2)\right]\\
 &+te^{-a_0x_1}E_{\alpha,2}(\gamma t^\alpha)\left[ \mu_1sin(\sqrt{b_0}x_2)+\mu_2cos(\sqrt{b_0}x_2)\right] \text{if}\ \alpha\in(1,2],
\end{aligned}
 \end{eqnarray}
where $\gamma=(a_0^2c_0-b_0\beta_0-a_0d_0+k_1)$ and $a_0,c_0,b_0,\beta_0,d_0,k_1,\mu_i, \nu_i \in \mathbb{R}$, $i=1,2$. Observe that the obtained exact solutions of \eqref{5.10} satisfy the given initial conditions \eqref{ic1} and \eqref{ic2}.

Now, we consider the another initial value problem for two-dimensional time-fractional convection-reaction-diffusion-wave equation \eqref{1} with cubic non-linearity in the form
\begin{eqnarray}
\begin{aligned}\label{5.12}
\frac{\partial^{\alpha}u}{\partial t^{\alpha}}= &c_0\frac{\partial^2 u }{\partial x_{1}^{2}}
+\frac{\partial}{\partial x_2}\left[(\beta_2u^2+\beta_1u+\beta_0)\frac{\partial u }{\partial x_{2}}\right]+(3b_0\beta_2u^2+\lambda_1u+\lambda_0)\frac{\partial u }{\partial x_2}\\
&+(b_0\lambda_1-2b_0^{2}\beta_1)u^2+k_1u,\ \alpha\in(0,2],\ t\geq0,
\end{aligned}
\end{eqnarray}
with the initial conditions
\begin{eqnarray}
&&\label{iv3} \text{If}\ \alpha \in(0,1], \text{then}\; u(x_1,x_2,0)=e^{-b_0x_2}\left(\nu_1sin(\sqrt{a_0}x_1)+\nu_2cos(\sqrt{a_0}x_1)\right).\\
&&\label{iv4} \text{If}\ \alpha \in(1,2], \text{then} \left\{
              \begin{array}{ll}
                u(x_1,x_2,0)=e^{-b_0x_2}\left(\nu_1sin(\sqrt{a_0}x_1)+\nu_2cos(\sqrt{a_0}x_1)\right)\ \& & \\
                {\dfrac{\partial u}{\partial t}}\mid_{t=0} =e^{-b_0x_2}\left(\mu_1sin(\sqrt{a_0}x_1)+\mu_2cos(\sqrt{a_0}x_1)\right). &
              \end{array}
            \right.
 \end{eqnarray}
We note that the above equation \eqref{5.12} admits a two-dimensional invariant subspace
 $\mathcal{V}_2=Span\{e^{-b_0x_2}sin(\sqrt{a_0}x_1),e^{-b_0x_2}cos(\sqrt{a_0}x_1)\}$ which was discussed in case 6 of table 10 (with $d_0=0$). Then by using the above-similar procedure, the obtained exact solutions for the two-dimensional time-fractional convection-reaction-diffusion-wave equation \eqref{5.12} are in the form
\begin{eqnarray}
\begin{aligned}\label{5.13}
 u(x_1,x_2,t)=&e^{-b_0x_2}E_{\alpha,1}(\gamma t^\alpha)\left(\nu_1sin(\sqrt{a_0}x_1)+\nu_2cos(\sqrt{a_0}x_1)\right)\;\text{if}\;\alpha \in (0,1], \\
 u(x_1,x_2,t)=&e^{-b_0x_2}E_{\alpha,1}(\gamma t^\alpha)\left(\nu_1sin(\sqrt{a_0}x_1)+\nu_2cos(\sqrt{a_0}x_1)\right)\\
 &+e^{-b_0x_2}E_{\alpha,2}(\gamma t^\alpha)\left(\mu_1sin(\sqrt{a_0}x_1)+\mu_2cos(\sqrt{a_0}x_1)\right)\;\text{if}\;\alpha\in (1,2],
\end{aligned}
\end{eqnarray}
where  $\gamma=(-a_0c_0+b_0^2\beta_0-b_0\lambda_0+k_1) $ and $a_0,b_0,c_0,\beta_0,\lambda_0, \nu_i, \mu_i \in \mathbb{R}$, $i=1,2$. Observe that the above solutions \eqref{5.13} satisfy the given initial conditions \eqref{iv3} and \eqref{iv4}.
\subsection{Exact solutions for the given equation \eqref{1} with quadratic non-linearity using the polynomial subspace}
 Now, let us consider the initial value problem for two-dimensional time-fractional convection-reaction-diffusion-wave equation \eqref{1} with quadratic non-linearity
\begin{eqnarray}
\begin{aligned}\label{5.20}
\frac{\partial^{\alpha}u}{\partial t^{\alpha}}=&\frac{\partial }{\partial x_1}\left[(c_1u+c_0)\frac{\partial u }{\partial x_1}\right]+\frac{\partial}{\partial x_2}\left[(\beta_1u+\beta_0)\frac{\partial u }{\partial x_{2}}\right] \\
&+d_0\frac{\partial u }{\partial x_1}+\lambda_0\frac{\partial u }{\partial x_2}+k_0,\ \alpha\in(0,2],\ t\geq0,
\end{aligned}
\end{eqnarray}
along with the initial conditions
\begin{eqnarray}
&&u(x_1,x_2,0)= \nu_1+\nu_2x_1+\nu_3x_2 \;\textsl{if}\; \alpha \in(0,1]\label{iv1},\\
&&\label{iv2}u(x_1,x_2,0)=\nu_1+\nu_2x_1+\nu_3x_2 \;\&\;{\dfrac{\partial u}{\partial t}}\mid_{t=0} =\mu_1+\mu_2x_1+\mu_3x_2 \; \textsl{if}\; \alpha \in(1,2].
\end{eqnarray}
 We know that the equation \eqref{5.20} admits a three-dimensional invariant subspace $\mathcal{V}_3=\text{Span}\{1,x_1,x_2\}$ which is discussed in case 3 of table 4.  Thus, the obtained exact solutions for the two-dimensional time-fractional convection-reaction-diffusion-wave equation \eqref{1} with quadratic non-linearity having the form
\begin{eqnarray}
\begin{aligned}
  u(x_1,x_2,t)=&\gamma_0\frac{t^\alpha}{\Gamma(\alpha+1)}+\nu_1+\nu_2x_1+\nu_3x_2  \;\text{if} \; \alpha \in (0,1],\\
  u(x_1,x_2,t)=&\gamma_0\frac{t^\alpha}{\Gamma(\alpha+1)}+\gamma_1\frac{t^{\alpha+1}}{\Gamma(\alpha+2)}+2\gamma_2\frac{t^{\alpha+2}}{\Gamma(\alpha+3)}+t\mu_1+\nu_1\\&+(\nu_2+t\mu_2)x_1+(\nu_3+t\mu_3)x_2\;\text{if} \;\alpha \in (1,2],
\end{aligned}
\end{eqnarray}
where $\gamma_0= c_1\nu_{2}^2+\beta_1\nu_3^2+d_0\nu_2+\lambda_0\nu_3+k_0,
\gamma_1=2(c_1\nu_2\mu_2+\beta_1\nu_3\mu_3)+d_0\mu_2+\lambda_0\mu_3,
\gamma_2=c_1\nu_2^{2}+\beta_1\nu_3^2$ and $ c_1,\beta_1,d_0,\lambda_0,k_0,c_0,\beta_0,\nu_i,\mu_i\in\mathbb{R}$, $i=1,2.$ It should be noted that the above solutions satisfy the initial conditions \eqref{iv1} and \eqref{iv2}.

\section{Invariant subspace method to two-dimensional time-fractional non-linear PDE with time delay}
In this section, we explain how to extend the invariant subspace method to two-dimensional time-fractional non-linear delay PDEs. In addition, we also explain how to derive the exact solution for the initial value problem of the two-dimensional time-fractional convection-reaction-diffusion-wave equation with linear term involving time delay.
\subsection{Invariant subspace method to two-dimensional time-fractional non-linear PDE with linear time delay}
In this subsection, we give a detailed study for constructing the invariant subspaces of the two-dimensional time-fractional non-linear PDEs with linear time delay. Thus, we consider the generalized two-dimensional time-fractional non-linear PDE with linear terms involving several time delay
\begin{eqnarray}\label{D1}
\begin{aligned}
& \frac{\partial^{\alpha}u}{\partial t^{\alpha}}=\tilde{\mathcal{K}}[u,\hat{u_i}]\equiv\mathcal{\hat{K}}[u] + \sum_{i=1}^{q}\mu_i\hat{u_i},\ \alpha>0,\ t>0,\\
& u(x_1,x_2,t)= \vartheta(x_1,x_2,t) \quad \text{if} \; t \in [-\hat{\tau},0],
\end{aligned}
\end{eqnarray}
{where} $u=u(x_1,x_2,t)$, $\hat{u_i}=u(x_1,x_2,t-\tau_i),\ \tau_i>0, \hat{\tau}=max \{ \tau_i: i=1,2,..,q \}$, $x_1,x_2\in\mathbb{R}$, $\mu_i>0, i=1,2,\dots,q$, $q\in \mathbb{N}$, $\dfrac{\partial^{\alpha}}{\partial t^{\alpha}}(\cdot) $ denotes the Caputo fractional derivative \eqref{c} of order $\alpha$, and $\mathcal{\hat{K}}[u]$ is the sufficiently given smooth differential operator of order $k$, that is,
\begin{equation}
 \mathcal{\hat{K}}[u]=\mathcal{\hat{K}}\left( x_{1},x_{2},
\frac{\partial u}{\partial x_1},\frac{\partial u}{\partial x_2},\frac{\partial^2 u}{\partial x_1^2},\frac{\partial^2 u}{\partial x_2^2},\frac{\partial^{2}u}{\partial x_1\partial x_2},\ldots,\frac{\partial^k u}{\partial x_1^k},\frac{\partial^k u}{\partial x_2^k},\frac{\partial^{k}u}{\partial x_1^{k_1}\partial x_2^{k_2}}\right)\label{D2}
\end{equation}
 and $k_1+k_2=k$, $k_1,k_2\in\mathbb{N}$.

Then the finite-dimensional linear space $\mathcal{V}_n$ given in \eqref{1.3} is said to be invariant under the differential operator $\mathcal{\tilde{K}}[u,\hat{u_i}]$ if for every $ u\in \mathcal{V}_n$ implies $\mathcal{\tilde{K}}[u,\hat{u}_i] \in \mathcal{V}_n$, which can be written as
 \begin{eqnarray}\begin{aligned}
  \mathcal{\tilde{K}}\left[\sum_{m=1}^{n}\kappa_m\xi_m(x_1,x_2), \sum_{m=1}^{n}\tilde{\kappa}_m\xi_m(x_1,x_2)\right]=&\sum \limits^{n} _{m=1}\Psi_m(\kappa_1,\kappa_2,\ldots,\kappa_n)\xi_m(x_1,x_2) \label{D3}\\
 & +\sum_{i=1}^{q}\sum_{m=1}^{n} \mu_i\tilde{\kappa}_m\xi_m(x_1,x_2),
\end{aligned}
 \end{eqnarray}
where $\tilde{\kappa}_m,\kappa_m \in \mathbb{R}$ and $\Psi_m $ denotes coefficient of expansion $\mathcal{\tilde{K}}[u,\hat{u}_i] $ with respect to the basis in \eqref{1.2}, $m=1,2,\dots,n$.
\begin{thm}
Suppose that the finite-dimensional linear space $\mathcal{V}_n$ defined in \eqref{1.3}, is invariant under the non-linear differential operator $\mathcal{\tilde{K}}[u,\hat{u_i}]$ given in the equation \eqref{D1}, then the two-dimensional non-linear time-fractional PDE with several linear time delays \eqref{D1} admits an exact solution of the form
\begin{equation}\label{D.4}
u(x_1,x_2,t)=\sum_{m=1}^{n}\Phi_m(t)\xi_m(x_1,x_2),
 \end{equation}
where the functions $\Phi_m(t)$ satisfy the system of ODEs of fractional-order
 \begin{equation}
\frac{d^\alpha \Phi_m(t)}{dt^\alpha}=\Psi_m( \Phi_1(t),\ldots,\Phi_n(t))+\sum_{i=1}^{q}\mu_i\Phi_m(t-\tau_i), \quad  m=1,2,\dots,n.
\label{D.5}
\end{equation}
\end{thm}
\textbf{Proof.}
Suppose that $\mathcal{V}_n$ be an n-dimensional invariant subspace admitted by the given differential operator ${\mathcal{\tilde{K}}}[u,\hat{u_i}].$
 Now, let us assume that $u(x_1,x_2,t)=\sum\limits_{m=1}^{n}\Phi_m(t)\xi_m(x_1,x_2)$.
Computing the Caputo fractional derivative of $u(x_1,x_2,t)$ of order $\alpha>0$ with respect to $t$ gives
 \begin{equation}\label{D.6}
\dfrac{\partial^{\alpha}u}{\partial t^{\alpha}}=\sum^{n} _{m=1}  \frac{d^{\alpha}\Phi_m(t)}{dt^\alpha}\xi_m(x_1,x_2).
 \end{equation}
Since $\mathcal{V}_n$ is invariant under ${\mathcal{\tilde{K}}}[u,\hat{u}_{i}]$. Thus, we have
\begin{eqnarray}\begin{aligned}
\mathcal{\tilde{K}}[u,\hat{u_i}]=&\mathcal{\tilde{K}}\left[u(x_1,x_2,t),{u}(x_1,x_2,t-\tau_i)\right]\\ =&\mathcal{\tilde{K}}\left[\sum \limits^{n} _{m=1}\Phi_m(t)\xi(x_1,x_2),\sum \limits^{n} _{m=1}\Phi_m(t-\tau_i)\xi(x_1,x_2,t)\right]\\
=&\sum \limits^{n} _{m=1}\Psi_m(\Phi_1(t),\Phi_2(t),...,\Phi_n(t))\xi_m(x_1,x_2)\\&+\sum \limits^{n} _{m=1}\sum \limits^{q} _{i=1}\mu_i{\Phi}_m(t-\tau_i)\xi_m(x_1,x_2).\label{D.7}
\end{aligned}
\end{eqnarray}
Substituting \eqref{D.6} and \eqref{D.7} in \eqref{D1}, we get
  \begin{equation}\label{D.8}
    \sum\limits_{m=1}^{n}\left[\frac{d^{\alpha}\Phi_m(t)}{dt^\alpha}-\Psi_m(\Phi_1(t),\Phi_2(t),...,\Phi_n(t))-\sum \limits^{q} _{i=1}\mu_i{\Phi}_m(t-\tau_i)\right]\xi_m(x_1,x_2) =0.
   \end{equation}
Since the functions $\xi_m(x_1,x_2)$, $i=1,\dots,m$ are linearly independent. Thus, the above equation \eqref{D.8} reduces to system of fractional ODEs \eqref{D.5}.

\textbf{Estimation of type I and type II invariant subspaces for \eqref{D1}:}
 Here we explain how to find the type I and type II linear spaces for the differential operator $\mathcal{\tilde{K}}[u,\hat{u}_i]$ that is given in \eqref{D1}. First, we consider the linear spaces which are defined as the solution space of homogeneous linear ODEs
\begin{eqnarray*}
\begin{aligned}
&\mathcal{V}_{n_1}=\left\{y_1\ \big{|}\ \mathcal{D}^{n_1}_{x_1}[y_1]\equiv \frac{d^{n_1}y_1}{dx_1^{n_1}} +a_{n_1-1}\frac{d^{n_1-1}y_1}{dx_1^{n_1-1}}+a_{n_1-2}\frac{d^{n_1-2}y_1}{dx_1^{n_1-2}}+\dots+a_{0}y_1=0\right\},\\
&\mathcal{V}_{n_2}=\left\{y_2\ \big{|}\ \mathcal{D}^{n_2}_{x_2}[y_2]\equiv \frac{d^{n_2}y_2}{dx_2^{n_2}} +b_{n_2-1}\frac{d^{n_2-1}y_2}{dx_2^{n_2-1}}+b_{n_2-2}\frac{d^{n_2-2}y_2}{dx_2^{n_2-2}}+\dots+b_{0}y_2=0\right\},
\end{aligned}
\end{eqnarray*}
where the functions $y_i$ are the linear combinations of the $n_i$-linearly independent solutions of $\mathcal{D}_{x_i}^{n_i}[y_i]=0$, $i=1,2$. Then, let us assume that the functions $v_1(x_1),v_2(x_2),\dots,v_{n_1}(x_1)$ are $n_1$-linearly independent solutions of $\mathcal{D}_{x_1}^{n_1}[y_1]=0$. Similarly, we assume that another set of functions $\zeta_1(x_2),\zeta_2(x_2),\dots,\zeta_{n_1}(x_2)$ are $n_2$-linearly independent solutions of $\mathcal{D}_{x_2}^{n_2}[y_2]=0$.
Then we can find two-types of linear spaces for the given differential operator $\mathcal{\tilde{K}}[u,\hat{u}_{i}]$. Suppose the given differential operator $\mathcal{\tilde{K}}[u,\hat{u}_{i}]$ admits type I and type II linear spaces $\mathcal{V}_{n_1n_2}$ and $\mathcal{V}_{n_1+n_2-1}$ that are given in \eqref{1.9}and \eqref{1.10}. Then the invariant conditions for the type I and type II  invariant subspaces are obtained as follows
\begin{itemize}
\item[1.] Type I invariance conditions for $\mathcal{K}[u,\hat{u_i}]$ are obtained in the following form
\begin{eqnarray*}
\begin{aligned}
&\mathcal{D}^{n_1}_{x_{1}}[\mathcal{\tilde{K}}]\equiv \frac{d^{n_1}\mathcal{\tilde{K}}}{dx_1^{n_1}}+a_{n_1-1}\frac{d^{n_1-1}\mathcal{\tilde{K}}}{dx_1^{n_1-1}} +\dots+a_1\frac{d\mathcal{\tilde{K}}}{dx_1}+a_{0}\mathcal{\tilde{K}}=0 \ \&\\
&\mathcal{D}^{n_2}_{x_2}[\mathcal{\tilde{K}}]\equiv \frac{d^{n_2}\mathcal{\tilde{K}}}{dx_2^{n_2}}+b_{n_2-1}\frac{d^{n_2}\mathcal{\tilde{K}}}{dx_2^{n_2}}+\dots+b_1\frac{d\mathcal{\tilde{K}}}{dx_2}+b_{0}\mathcal{\tilde{K}}=0 \label{t1}\\
&\text{along with}\\
&\mathcal{D}^{n_1}_{x_{1}}[u]\equiv \frac{d^{n_1}u}{dx_1^{n_1}}+a_{n_1-1}\frac{d^{n_1-1}u}{dx_1^{n_1-1}}+\dots+a_1\frac{du}{dx_1}+a_{0}u=0\; \&\;\\
&\mathcal{D}^{n_2}_{x_2}[u]\equiv \frac{d^{n_2}u}{dx_2^{n_2}}+b_{n_2-1}\frac{d^{n_2-1}u}{dx_2^{n_2-1}}+\dots+b_1\frac{du}{dx_2}+b_{0}u=0,
\end{aligned}
\end{eqnarray*}
 where $a_i, \ (i=0,1,\dots,n_{1}-1)$ and $ b_j, \  (j=0,1,\dots,n_{2}-1)$ are constants to be determined.
\item[2] Type II invariance conditions for $\mathcal{K}[u,\hat{u_i}]$ take the form
\begin{eqnarray*}\begin{aligned}\label{t2}
&\mathcal{D}^{n_1}_{x_{1}}[\mathcal{\tilde{K}}]\equiv \frac{d^{n_1}\mathcal{\tilde{K}}}{dx_1^{n_1}}+a_{n_1-1}\frac{d^{n_1-1}\mathcal{\tilde{K}}}{dx_1^{n_1-1}}+\dots+a_1\frac{d\mathcal{\tilde{K}}}{dx_1}=0,\\
&\mathcal{D}^{n_2}_{x_2}[\mathcal{\tilde{K}}]\equiv \frac{d^{n_2}\mathcal{\tilde{K}}}{dx_2^{n_2}}+b_{n_2-1}\frac{d^{n_2-1}\mathcal{\tilde{K}}}{dx_2^{n_2-1}}+\dots+b_1\frac{d\mathcal{\tilde{K}}}{dx_2}=0\ \& \ \frac{\partial^{2}\mathcal{\tilde{K}}}{\partial x_1\partial x_2}=0\;\\
&\text{along with}\\
&\mathcal{D}^{n_1}_{x_{1}}[u]\equiv \frac{d^{n_1}u}{dx_1^{n_1}}+a_{n_1-1}\frac{d^{n_1-1}u}{dx_1^{n_1-1}}+\dots+a_1\frac{du}{dx_1}=0,\\
&\mathcal{D}^{n_2}_{x_2}[u]\equiv \frac{d^{n_2}u}{dx_2^{n_2}}+b_{n_2-1}\frac{d^{n_2-1}u}{dx_2^{n_2-1}}+\dots+b_1\frac{du}{dx_2}=0\ \&\ \frac{\partial^{2} u(x_1,x_2,t)}{\partial x_1\partial x_2}=0,
\end{aligned}
\end{eqnarray*}
 where $a_i, \ (i=1,\dots,n_{1}-1)$ and  $b_j,\ (j=1,\dots,n_{2}-1)$ are constants to be determined.\\
 Next, the applicability and effectiveness of the method have been illustrated through the two-dimensional time-fractional cubic non-linear convection-reaction-diffusion-wave equation with linear time delay.
\end{itemize}
\subsection{Exact solution for the two-dimensional time-fractional cubic non-linear convection-reaction-diffusion-wave equation with linear time delay}
In this subsection, let us consider the two-dimensional time-fractional cubic non-linear convection-reaction-diffusion-wave equation with linear time delay
\begin{eqnarray}
\begin{aligned}\label{6.9}
\frac{\partial^{\alpha}u}{\partial t^{\alpha}}=&\frac{\partial }{\partial x_1}\left[\left(\frac{d_2}{3a_0}u^2+c_1u+c_0\right)\frac{\partial u }{\partial x_1}\right]+\beta_0\frac{\partial^2 u }{\partial x_{2}^{2}}+(d_2u^2+d_1u+d_0)\frac{\partial u }{\partial x_1}\\
&+\lambda_0\frac{\partial u }{\partial x_2}+(-2a_0^2c_1+a_0d_1)u^2+k_1u+\mu\hat{u},\ \alpha\in(0,2],\ t\geq0,\\
u(x_1,x_2,t)=&\vartheta(x_1,x_2,t)=\phi_1(t)e^{-a_0x_1}+\phi_2(t)e^{-(a_0x_1+b_1x_2)},\;t \in [-\tau,0],
\end{aligned}
\end{eqnarray}
along with the initial conditions
\begin{eqnarray}
&&\label{ii1}u(x_1,x_2,0)=\kappa_1e^{-a_0x_1}+\kappa_2 e^{-(a_0x_1+b_1x_2)} \;  \textsl{if}\; \alpha \in(0,1] ,\; t \in [-\tau,0],\\
&&\nonumber u(x_1,x_2,0)=\kappa_1e^{-a_0x_1}+\kappa_2e^{-(a_0x_1+b_1x_2)}\ \& \\
&&\label{ii2}\frac{\partial u}{\partial t}\mid_{t=0}=\hat{\kappa}_1e^{-a_0x_1}+\hat{\kappa}_2e^{-(a_0x_1+b_1x_2)} \;\textsl{if}\; \alpha \in(1,2,]\label{ic}
\end{eqnarray}
where $\hat{u}=u(x_1,x_2,t-\tau), \tau>0,\kappa_i,\hat{\kappa}_i,\mu\in\mathbb{R},\ i=1,2$.\\
Here, we consider the differential operator $\mathcal{\hat{K}}[u,\hat{u}]$ for the given equation \eqref{6.9} as
\begin{eqnarray}
\begin{aligned}
\mathcal{\hat{K}}[u,\hat{u}]=&\frac{\partial }{\partial x_1}\left[\left(\frac{d_2}{3a_0}u^2+c_1u+c_0\right)\frac{\partial u }{\partial x_1}\right]+\beta_0\frac{\partial^2 u }{\partial x_{2}^{2}}+(d_2u^2+d_1u+d_0)\frac{\partial u }{\partial x_1}\\&+\lambda_0\frac{\partial u }{\partial x_2}+(-2a_0^2c_1+a_0d_1)u^2+k_1u+\mu\hat{u}.
\end{aligned}
\end{eqnarray}
It should be noted that the cubic non-linear differential operator $\mathcal{\hat{K}}[u,\hat{u}]$ admits a two-dimensional exponential linear space $\mathcal{V}_2=Span\{e^{-a_0x_1},e^{-(a_0x_1+b_1x_2)}\}$, since for some constants $\delta_k,\hat{\delta}_k\in \mathbb{R},k=1,2$
\begin{eqnarray*}\begin{aligned}
&\mathcal{\hat{K}}[\delta_1e^{-a_0x_1}+\delta_2e^{-(a_0x_1+b_1x_2)},\hat{\delta}_1e^{-a_0x_1}+\hat{\delta}_2e^{-(a_0x_1+b_1x_2)}]\\=&(a_0^2c_0-d_0a_0+k_1)\delta_1e^{-a_0x_1}+
(k_1+\beta_0b_1^2-\lambda_0b_1+c_0a_0^2-a_{0}d_{0})\delta_2e^{-(a_0x_1+b_1x_2)}\\
&+\mu[\hat{\delta}_1e^{-a_0x_1}+\hat{\delta}_2e^{-(a_0x_1+b_1x_2)}]\in\mathcal{V}_2,
\end{aligned}
\end{eqnarray*}
 which suggests an exact solution of \eqref{6.9} in the form
\begin{equation}\label{6.10}
u(x_1,x_2,t)=\Phi_1(t)e^{-a_0x_1}+\Phi_2(t)e^{-(a_0x_1+b_1x_2)},
\end{equation}
where the functions $\Phi_1(t)$ and $\Phi_2(t)$ satisfy the following system
\begin{equation}\label{6.11}
\frac{d^{\alpha}\Phi_i(t)}{d t^{\alpha}}=\gamma_i{\Phi_i(t)}+\mu\Phi_i(t-\tau),\ i=1,2,\ \alpha\in(0,2],
\end{equation}
where $\gamma_1=(a_0^2c_0-d_0a_0+k_1) $ and $\gamma_2=(k_1+\beta_0b_1^2-\lambda_0b_1+c_0a_0^2-a_{0}d_{0})$.
Now, let us first consider the case $\alpha\in(0,1]$. Applying the Laplace transformation of equation \eqref{6.11}, that reads as
\begin{equation*}
\mathcal{L}\left[\frac{d^{\alpha}\Phi_{i}(t)}{d t^{\alpha}}\right]=\mathcal{L}\left[\gamma_{i}{\Phi_{i}(t)}\right]+\mu \mathcal{L}[\Phi_i(t-\tau)], i=1,2.
\end{equation*}
 From the given initial data \eqref{6.9}, we get $ \Phi_i(t)=\phi_i(t),\forall t \in [-\tau,0],i=1,2.$ Thus, we obtain
\begin{eqnarray*}
s^\alpha\tilde{\Phi}_i(s)-s^{\alpha-1}\Phi_i(0)=\gamma_i\tilde{\Phi}_i(s)+\mu e^{-\tau s}\tilde{\Phi}_i(s)+\mu e^{-\tau s}{\int_{-\tau}^{0}e^{-s\xi}\phi_i(\xi)d\xi}.
\end{eqnarray*}
The above equation can be written as
\begin{equation}\label{iv6}
\tilde{\Phi}_i(s)=\Phi_i(0)\left(\frac{s^{\alpha-1}}{s^\alpha-\gamma_i-\mu e^{-\tau s}}\right)+\left(\frac{\mu e^{-\tau s}}{s^\alpha-\gamma_i-\mu e^{-\tau s}}\right){\int_{-\tau}^{0}e^{-s\xi}\phi_i(\xi)d\xi}.
\end{equation}
Taking the inverse Laplace transformation of \eqref{iv6}, we get
\begin{eqnarray}\label{iv5}
\Phi_i(t)=\Phi_i(0)\mathcal{ L}^{-1}\left[\frac{s^{\alpha-1}}{s^\alpha-\gamma_i-\mu e^{-\tau s}}\right]+\mu\mathcal{ L}^{-1}\left[\left(\frac{e^{-\tau s}}{s^\alpha-\gamma_i-\mu e^{-\tau s}}\right){\int_{-\tau}^{0}e^{-s\xi}\phi_i(\xi)d\xi}\right].
\end{eqnarray}
Let $0<\mid{\dfrac{\mu e^{-\tau s}}{s^\alpha-\gamma_{i}}}\mid<1$. Now, we can simplify the first term of \eqref{iv5} as
\begin{eqnarray*}
\begin{aligned}
\mathcal{ L}^{-1}\left[\frac{s^{\alpha-1}}{s^\alpha-\gamma_i-\mu e^{-\tau s}}\right]=&\mathcal{ L}^{-1}\left[\sum_{m=0}^{\infty}\mu^me^{-\tau m s}\frac{s^{\alpha-1}}{(s^\alpha-\gamma_i)^{m+1}}\right] \\
 =&\sum_{m=0}^{\infty}\mu^m\mathcal{U}_{m\tau}(t)(t-m\tau)^{\alpha m} E_{\alpha,\alpha m+1}^{m+1}(\gamma_i(t-m\tau)^{\alpha}),
 \end{aligned}
\end{eqnarray*}
where  $\mathcal{U}_a(t)$ denotes the unit step function that can be defined as
 $\mathcal{U}_a(t)=\left\{
         \begin{array}{ll}

         1 \quad t\geq a\\
       0\quad t< a
         \end{array}
       \right.$
       \\
 and $E_{p,q}^{r}(\cdot)$ is the generalized three-parameter Mittag-Leffler function \cite{mathai}, is defined as $ E_{p,q}^{r}(w)=\sum\limits_{m=0}^{\infty}\dfrac{(r)_m w^m}{\Gamma(pm+q) m!}$, $ w\in \mathbb{R}$, $(r)_m=\dfrac{\Gamma(r+m)}{\Gamma(r)}$, $p,q,r>0$.
\\ Next, by using the Laplace convolution theorem which states that any two piecewise continuous functions $ f_1(t)$ and  $f_2(t)$ defined on $[0,\infty)$ and of exponential order $v>0$, then Laplace transformation of convolution of $f_1(t)$ and $f_2(t)$ is given by
\[ \mathcal{L}\left[ (f_1\ast f_2)(t)\right]=\mathcal{L}(f_1(t))\mathcal{L}(f_2(t)), \] where $ (f_1\ast f_2)(t)=\int_{0}^{t}{f_1(\xi)f_2(t-\xi)d\xi}=\int_{0}^{t}{f_1(t-\xi)f_2(\xi)d\xi}.$ Also,
 consider the extension of $ \phi_i(t) $ on $[-\tau,\infty]$ of the form $ \phi_i(t)=\left\{
         \begin{array}{ll}
                  \Phi_i(t) \quad t\in [-\tau,0] \\
       \Phi_i(0) \quad t\geq 0
        \end{array}
      \right.,i=1,2$ and define the function $f(t)=\left\{
         \begin{array}{ll}
                 0 \quad t\geq 0\\
       1\quad t< 0
         \end{array}
       \right.$,
then we can simplify the second term\\ $\mathcal{ L}^{-1}\left[\left(\dfrac{e^{-\tau s}}{s^\alpha-\gamma_i-\mu e^{-\tau s}}\right){\int_{-\tau}^{0}e^{-s\xi}\phi_i(\xi)d\xi}\right]$ of \eqref{iv5} for $0<\mid{\dfrac{\mu e^{-\tau s}}{s^\alpha-\gamma_{i}}}\mid<1$, $i=1,2$, as follows
\begin{eqnarray*}
\begin{aligned}
 &\mathcal{ L}^{-1}\left[\left(\frac{e^{-\tau s}}{s^\alpha-\gamma_i-\mu e^{-\tau s}}\right){\int_{-\tau}^{0}e^{-s\xi}\phi_i(\xi)d\xi}\right]=\sum_{m=0}^{\infty}\mathcal{ L}^{-1}\left[\frac{\mu^m e^{-\tau(m+1) s}}{(s^\alpha-\gamma_i)^{m+1}}{\int_{-\tau}^{0}e^{-s\xi}\phi_i(\xi)d\xi}\right]\\
=&\sum_{m=0}^{\infty}\mu^m \mathcal{ L}^{-1}\left[\frac{e^{-\tau m s}}{(s^\alpha-\gamma_i)^{m+1}}\right]\ast\mathcal{ L}^{-1}\left[e^{-\tau  s}{\int_{-\tau}^{0}e^{-s\xi}\phi_i(\xi)d\xi}\right]\\
=&\sum_{m=0}^{\infty}\mu^m\mathcal{ L}^{-1}\left[\frac{e^{-\tau m s}}{(s^\alpha-\gamma_i)^{m+1}}\right]\ast\mathcal{ L}^{-1}\left[\int_{0}^{\infty}e^{-s\xi}\phi_i(\xi-\tau)f(\xi-\tau)d\xi\right]\\
=&\sum_{m=0}^{\infty}\mu^m\mathcal{U}_{m\tau}(t)(t-m\tau)^{\alpha( m+1)-1} E_{\alpha,\alpha( m+1)}^{m+1}(\gamma_i(t-m\tau)^{\alpha})\ast[\phi_i(t-\tau)f(t-\tau)].
\end{aligned}
\end{eqnarray*}
Thus, we obtain the function $\Phi_i(t)$, $i=1,2$, in the form
\begin{eqnarray*}\begin{aligned}
\Phi_i(t)=&\Phi_i(0)\sum\limits_{m=0}^{\infty}\mu^m\mathcal{U}_{m\tau}(t)(t-m\tau)^{\alpha m} E_{\alpha,\alpha m+1}^{m+1}(\gamma_i(t-m\tau)^{\alpha})\\ +&
  \left[ \sum\limits_{m=0}^{\infty}\mu^{m+1}\mathcal{U}_{m\tau}(t)(t-m\tau)^{\alpha( m+1)-1} E_{\alpha,\alpha( m+1)}^{m+1}(\gamma_i(t-m\tau)^{\alpha})\right]\ast[\phi_i(t-\tau)f(t-\tau)].\end{aligned}
\end{eqnarray*}
Hence the obtained exact solution of the two-dimensional time-fractional cubic non-linear convection-reaction-diffusion equation with linear time delay \eqref{6.9} for
$0<\mid{\dfrac{\mu e^{-\tau s}}{s^\alpha-\gamma_{i}}}\mid<1,i=1,2,$ as follows
\begin{equation}
\begin{split}\label{ii3}
 u(x_1,x_2,t)&=e^{-a_0x_1}\left\{\kappa_1\sum_{m=0}^{n}\mu^m(t-m\tau)^{\alpha m} E_{\alpha,\alpha m+1}^{m+1}(\gamma_1(t-m\tau)^{\alpha}) \right.\\
 +&\left.       \left[\sum_{m=0}^{n}\mu^{m+1}(t-m\tau)^{\alpha( m+1)-1} E_{\alpha,\alpha( m+1)}^{m+1}(\gamma_1(t-m\tau)^{\alpha})\right]\ast[\phi_1(t-\tau)f(t-\tau)] \right\}\\
 +&e^{-(a_0x_1+b_1x_2)}\left\{ \kappa_2 \sum_{m=0}^{n}\mu^m(t-m\tau)^{\alpha m} E_{\alpha,\alpha m+1}^{m+1}(\gamma_2(t-m\tau)^{\alpha}) \right. \\
 +&\left.    \left[     \sum_{m=0}^{n}\mu^{m+1}(t-m\tau)^{\alpha( m+1)-1} E_{\alpha,\alpha( m+1)}^{m+1}(\gamma_2(t-m\tau)^{\alpha})\right]\ast[\phi_2(t-\tau)f(t-\tau)] \right\},
\end{split}
\end{equation}
where $\alpha\in(0,1]$, $n-1< \dfrac{t}{\tau}\leq n$, $t>0,\ \tau>0$, $\kappa_i=\Phi_i(0),i=1,2$, $n\in \mathbb{N}$, $f(t)=\left\{
         \begin{array}{ll}
                   0 \quad t\geq 0\\
       1\quad t< 0
         \end{array}
       \right.,$ $\phi_i(t)=\Phi_i(t),\ \forall t \in [-\tau,0],i=1,2,$ $\gamma_1=(a_0^2c_0-d_0a_0+k_1)$  and
       $\gamma_2=(k_1+\beta_0b_1^2-\lambda_0b_1+c_0a_0^2-a_{0}d_{0})$.
       \\
It is important to note that the above solution \eqref{ii3} satisfies the given initial condition \eqref{ii1}. The obtained solution \eqref{ii3} can be viewed as
       \begin{equation}
\begin{split}\label{vic1}
 u(x_1,x_2,t)&=e^{-a_0x_1}\left\{\kappa_1\left(E_{\alpha,1}^{1}(\gamma_1t^\alpha)+\sum_{m=1}^{n}\mu^m(t-m\tau)^{\alpha m} E_{\alpha,\alpha m+1}^{m+1}(\gamma_1(t-m\tau)^{\alpha})\right) \right.\\
 +&\left.      \int_{0}^{t}{\sum_{m=0}^{n}\mu^{m+1}(\xi-m\tau)^{\alpha( m+1)-1} E_{\alpha,\alpha( m+1)}^{m+1}(\gamma_1(\xi-m\tau)^{\alpha})\phi_1(t-\tau-\xi)f(t-\tau-\xi)d\xi} \right\}\\
 +&e^{-(a_0x_1+b_1x_2)}\left\{ \kappa_2\left(E_{\alpha,1}^{1}(\gamma_2t^\alpha)+\sum_{m=1}^{n}\mu^m(t-m\tau)^{\alpha m} E_{\alpha,\alpha m+1}^{m+1}(\gamma_2(t-m\tau)^{\alpha})\right) \right. \\
 +&\left.   \int_{0}^{t}{    \sum_{m=0}^{n}\mu^{m+1}(\xi-m\tau)^{\alpha( m+1)-1} E_{\alpha,\alpha( m+1)}^{m+1}(\gamma_2(\xi-m\tau)^{\alpha})\phi_2(t-\tau-\xi)f(t-\tau-\xi)d\xi} \right\}.
\end{split}
\end{equation}
Substutite $t=0$ in \eqref{vic1}, we have
\[ u(x_1,x_2,0)=\kappa_1e^{-a_0x_1}+\kappa_2 e^{-(a_0x_1+b_1x_2)} ,\]
since $n=0$ and $E_{\alpha,1}^{1}(\gamma_it^\alpha)=1\;\text{at}\;t=0,i=1,2 $.

Now we consider the second case $\alpha\in (1,2]$, proceeding in the similar way as explained above for $\alpha\in (0,1]$. For this case, we applying the Laplace transformation to equation \eqref{6.11} which yields
\[\mathcal{L}\left[\frac{d^{\alpha}\Phi_{i}(t)}{d t^{\alpha}}\right]=\mathcal{L}\left[\gamma_{i}{\Phi_{i}(t)}\right]+\mu \mathcal{L}[\Phi_i(t-\tau)],i =1,2,\]
The above equation can be simplified as
\[ s^\alpha\tilde{\Phi}_i(s)-s^{\alpha-1}\Phi_i(0)-s^{\alpha-2}\Phi_i'(0)=\gamma_i\tilde{\Phi}_i(s)+\mu e^{-\tau s}\tilde{\Phi}_i(s)+\mu e^{-\tau s}{\int_{-\tau}^{0}e^{-s\xi}\Phi_i(\xi)d\xi},i =1,2,\]
which can be written in the form
\begin{eqnarray*}
\begin{aligned}
 \tilde{\Phi}_i(s)=& \Phi_i(0)\left(\frac{s^{\alpha-1}}{s^\alpha-\gamma_i-\mu e^{-\tau s}}\right)+\Phi_i'(0)\left(\frac{s^{\alpha-2}}{s^\alpha-\gamma_i-\mu e^{-\tau s}}\right)\\&+\left(\frac{\mu}{s^\alpha-\gamma_i-\mu e^{-\tau s}}\right)e^{-\tau s}{\int_{-\tau}^{0}e^{-s\xi}\phi_i(\xi)d\xi}.
\end{aligned}\end{eqnarray*}
Thus, the obtained function $\Phi_i(t)$, for $0<\mid{\dfrac{\mu e^{-\tau s}}{s^\alpha-\gamma_{i}}}\mid<1,i=1,2$,  is of the form
  \begin{eqnarray*}\begin{aligned}
 \Phi_i(t)=&\Phi_i(0)\sum\limits_{m=0}^{\infty}\mu^m\mathcal{U}_{m\tau}(t)(t-m\tau)^{\alpha m} E_{\alpha,\alpha m+1}^{m+1}(\gamma_i(t-m\tau)^{\alpha}) \\
 & +\Phi_i'(0)\sum\limits_{m=0}^{\infty}\mu^m\mathcal{U}_{m\tau}(t)(t-m\tau)^{\alpha m+1} E_{\alpha,\alpha m+2}^{m+1}(\gamma_i(t-m\tau)^{\alpha}) \\&+
   \left[ \sum\limits_{m=0}^{\infty}\mu^{m+1}\mathcal{U}_{m\tau}(t)(t-m\tau)^{\alpha( m+1)-1} E_{\alpha,\alpha( m+1)}^{m+1}(\gamma_i(t-m\tau)^{\alpha})\right]\ast[\phi_i(t-\tau)f(t-\tau)].\end{aligned}
  \end{eqnarray*}
Hence the exact solution of the two-dimensional time-fractional cubic non-linear convection-reaction-diffusion equation with linear time delay \eqref{6.9}, for $0<\mid{\dfrac{\mu e^{-\tau s}}{s^\alpha-\gamma_{i}}}\mid<1,i=1,2$, as follows
\begin{equation}
\begin{split}\label{ii4}
 u(x_1,x_2,t)=&e^{-a_0x_1}\left\{ \kappa_1\sum_{m=0}^{n}\mu^m(t-m\tau)^{\alpha m} E_{\alpha,\alpha m+1}^{m+1}(\gamma_1(t-m\tau)^{\alpha}) \right.\\&+\hat{\kappa}_1 \sum_{m=0}^{n}\mu^m(t-m\tau)^{\alpha m+1} E_{\alpha,\alpha m+2}^{m+1}(\gamma_1(t-m\tau)^{\alpha})
  \\&\left.    +   \left[  \sum_{m=0}^{n}\mu^{m+1}(t-m\tau)^{\alpha( m+1)-1} E_{\alpha,\alpha( m+1)}^{m+1}(\gamma_1(t-m\tau)^{\alpha})\right]\ast[\phi_1(t-\tau)f(t-\tau)] \right\}
 \\&
+ e^{-(a_0x_1+b_1x_2)}\left\{ \kappa_2 \sum_{m=0}^{n}\mu^m(t-m\tau)^{\alpha m} E_{\alpha,\alpha m+1}^{m+1}(\gamma_2(t-m\tau)^{\alpha}) \right.\\&  + \hat{\kappa_2}\sum_{m=0}^{n}\mu^m(t-m\tau)^{\alpha m+1} E_{\alpha,\alpha m+2}^{m+1}(\gamma_2(t-m\tau)^{\alpha})
 \\&\left.  +  \left[     \sum_{m=0}^{n}\mu^{m+1}(t-m\tau)^{\alpha( m+1)-1} E_{\alpha,\alpha( m+1)}^{m+1}(\gamma_2(t-m\tau)^{\alpha})\right]\ast[\phi_2(t-\tau)f(t-\tau)] \right\},
\end{split}
\end{equation}
where  $ \phi_i(t)=\Phi_i(t),\ \forall t\in[-\tau,0],\kappa_i=\Phi_i(0), \hat{\kappa}_i=\Phi'_i(0),\  i=1,2, $
$ n-1 <\dfrac{t}{\tau}\leq n, n\in\mathbb{N}$, $\gamma_1=(a_0^2c_0-d_0a_0+k_1)\; \text{and} \;\gamma_2=(k_1+\beta_0b_1^2-\lambda_0b_1+c_0a_0^2-a_{0}d_{0})$. Additionally, we observe that the obtained exact solution \eqref{ii4} satisfies the given initial conditions \eqref{ii2}.
The above solution \eqref{ii4} can be written as
\begin{equation}
\begin{split}\label{vic2}
 u(x_1,x_2,t)&=e^{-a_0x_1}\left\{\kappa_1\left(E_{\alpha,1}^{1}(\gamma_1t^\alpha)+\sum_{m=1}^{n}\mu^m(t-m\tau)^{\alpha m} E_{\alpha,\alpha m+1}^{m+1}(\gamma_1(t-m\tau)^{\alpha})\right) \right.
 \\+&  \hat{\kappa}_1\left( tE_{\alpha,2}^{1}(\gamma_1t^{\alpha})+\sum_{m=1}^{n}\mu^m(t-m\tau)^{\alpha m+1} E_{\alpha,\alpha m+2}^{m+1}(\gamma_1(t-m\tau)^{\alpha})\right)
 \\
 +&\left.      \int_{0}^{t}{\sum_{m=0}^{n}\mu^{m+1}(\xi-m\tau)^{\alpha( m+1)-1} E_{\alpha,\alpha( m+1)}^{m+1}(\gamma_1(\xi-m\tau)^{\alpha})\phi_1(t-\tau-\xi)f(t-\tau-\xi)d\xi} \right\}
 \\
 +&e^{-(a_0x_1+b_1x_2)}\left\{ \kappa_2\left(E_{\alpha,1}^{1}(\gamma_2t^\alpha)+\sum_{m=1}^{n}\mu^m(t-m\tau)^{\alpha m} E_{\alpha,\alpha m+1}^{m+1}(\gamma_2(t-m\tau)^{\alpha})\right) \right.
 \\+&
   \hat{\kappa_2}\left(tE_{\alpha,2}^{1}(\gamma_2t^{\alpha})+\sum_{m=1}^{n}\mu^m(t-m\tau)^{\alpha m+1} E_{\alpha,\alpha m+2}^{m+1}(\gamma_2(t-m\tau)^{\alpha})\right)
  \\
 +&\left.   \int_{0}^{t}{    \sum_{m=0}^{n}\mu^{m+1}(\xi-m\tau)^{\alpha( m+1)-1} E_{\alpha,\alpha( m+1)}^{m+1}(\gamma_2(\xi-m\tau)^{\alpha})\phi_2(t-\tau-\xi)f(t-\tau-\xi)d\xi} \right\},
\end{split}
\end{equation}
which satisfies the initial condition at $t=0$
\[ u(x_1,x_2,0)=\kappa_1e^{-a_0x_1}+\kappa_2 e^{-(a_0x_1+b_1x_2)} ,\]
since $n=0$ and $E_{\alpha,1}^{1}(\gamma_it^\alpha)=1\;\text{at}\;t=0,i=1,2 $.
And also, the solution \eqref{vic2} satisfies second initial condition
\[\frac{\partial u}{\partial t}\mid_{t=0}=\hat{\kappa}_1e^{-a_0x_1}+\hat{\kappa}_2e^{-(a_0x_1+b_1x_2)},\] because
 \begin{eqnarray*}
&\dfrac{d}{dt}\left[ \int_{0}^{t}{\sum_{m=0}^{n}\mu^{m+1}(\xi-m\tau)^{\alpha( m+1)-1} E_{\alpha,\alpha( m+1)}^{m+1}(\gamma_i(\xi-m\tau)^{\alpha})\phi_i(t-\tau-\xi)f(t-\tau-\xi)d\xi}  \right]\big{|}_{t=0}\\&= {\sum_{m=0}^{n}\mu^{m+1}(t-m\tau)^{\alpha( m+1)-1} E_{\alpha,\alpha( m+1)}^{m+1}(\gamma_i(t-m\tau)^{\alpha})\phi_i(-\tau)f(-\tau)}\big{|}_{t=0}=0,
\end{eqnarray*}
 $\dfrac{d}{dt}\left[E_{\alpha,1}^{1}(\gamma_it^\alpha)\right]\big{|}_{t=0}=0$ and $ \dfrac{d}{dt}\left[tE_{\alpha,2}^{1}(\gamma_it^{\alpha})\right]\big{|}_{t=0}=1,i=1,2.$

Now, we explain how to extend the invariant subspace method to generalized two-dimensional time-fractional non-linear PDE with time delay.
\subsection{Extension of the invariant subspaces to generalized two-dimensional time-fractional non-linear PDE with time delay}
Consider the generalized two-dimensional time-fractional time delay non-linear PDE in the form
\begin{eqnarray}\begin{aligned}
& \frac{\partial^{\alpha}u}{\partial t^{\alpha}}=\tilde{\mathcal{H}}[u,\hat{u}], \alpha>0,\ t>0\\
& u(x_1,x_2,t)=\omega(x_1,x_2,t)\; \text{if}\; t\in[-\tau,0],\ \tau>0,\label{gen.dly}
\end{aligned}
\end{eqnarray}
where $u=u(x_1,x_2,t)$, $\hat{u}=u(x_1,x_2,t-\tau),$  $x_1,x_2\in\mathbb{R}$, $\dfrac{\partial^{\alpha}}{\partial t^{\alpha}}(\cdot) $ denotes Caputo fractional derivative \eqref{c} of order $\alpha$, and $\mathcal{\tilde{H}}[u,\hat{u}]$ is the sufficiently given smooth differential operator of order $k$, that is,
\begin{eqnarray}\begin{aligned}
\mathcal{\tilde{H}}[u,\hat{u}]=&\mathcal{\tilde{H}}\left( x_{1},x_{2},\label{gen.dlyop}\frac{\partial u}{\partial x_1},
\frac{\partial \hat{u}}{\partial x_1},
\frac{\partial u}{\partial x_2},
\frac{\partial \hat{u}}{\partial x_2},
\frac{\partial^2 u}{\partial x_1^2},
\frac{\partial^2 \hat{u}}{\partial x_1^2},
\frac{\partial^2 u}{\partial x_2^2},
\frac{\partial^2\hat{ u}}{\partial x_2^2},
\frac{\partial^{2}u}{\partial x_1\partial x_2},\ldots,\frac{\partial^k u}{\partial x_1^k},\right.\\
&\left.\quad
\frac{\partial^{2}\hat{u}}{\partial x_1\partial x_2},\ldots,\frac{\partial^k\hat{u}}{\partial x_1^k},
\frac{\partial^k u}{\partial x_2^k},\frac{\partial^k \hat{u}}{\partial x_2^k},
\frac{\partial^{k} u}{\partial x_1^{k_1}\partial x_2^{k_2}},\frac{\partial^{k}\hat{u}}{\partial x_1^{k_1}\partial x_2^{k_2}}
\right)
\end{aligned}
\end{eqnarray}
and  $k_1+k_2=k$, $k_1,k_2  \in \mathbb{N}$.
Suppose the given differential operator \eqref{gen.dlyop} admits the invariant subspace $\mathcal{V}_n $ given in \eqref{1.4}. Then there exists $n$ functions $ \Psi_m \;  (m=1,2,\dots,n)$ such that
\begin{eqnarray}
\begin{aligned}
\mathcal{\tilde{H}}[u,\hat{u}]=&\mathcal{\tilde{H}}\left[\sum_{m=1}^{n}\kappa_m\xi_m(x_1,x_2), \sum_{m=1}^{n}\hat{\kappa}_m\xi_m(x_1,x_2)\right]\\
=&\sum \limits^{n} _{m=1}\Psi_m(\kappa_1,\ldots,\kappa_n,\hat{\kappa}_1,\ldots,\hat{\kappa}_m)\xi_m(x_1,x_2),
\end{aligned}
\end{eqnarray}
where $\kappa_m ,\hat{\kappa}_m \in \mathbb{R},\  m=1,2,\dots,n$.
 \begin{thm}
Suppose that the linear space $\mathcal{V}_n$ given in \eqref{1.3} is  invariant under the non-linear differential operator $\mathcal{\tilde{H}}[u,\hat{u}]$  given in \eqref{gen.dlyop}, then the generalized two-dimensional non-linear time-fractional PDE with time delay \eqref{gen.dly} has admits an exact solution in the form
\begin{equation}
u(x_1,x_2,t)=\sum_{m=1}^{n}\Phi_m(t)\xi_m(x_1,x_2),\label{ds1}
 \end{equation}
where the functions $\Phi_m(t)$ satisfy the system of ODEs of fractional-order
 \begin{equation}
\dfrac{d^\alpha \Phi_m(t)}{dt^\alpha}=\Psi_m( \Phi_1(t),\ldots,\Phi_n(t),\Phi_1(t-\tau),\ldots,\Phi_n(t-\tau)), m=1,2,\dots,n.
\label{sys1}
\end{equation}
\end{thm}
\textbf{Proof.} Proof of the theorem is similar to above-theorem 5.1.

\section{Conclusion}
The presented work was investigated how we can extend the invariant subspace method to two-dimensional time-fractional non-linear PDEs. More precisely, the systematic study was given for constructing the various dimensions of the invariant subspaces for the two-dimensional time-fractional generalized convection-reaction-diffusion-wave equation along with the initial conditions \eqref{1} for the first time. Also, we have shown explicitly that the time-fractional convection-reaction-diffusion-wave equation attains more than one invariant subspaces in the same dimension of linear spaces under consideration. Additionally, the special types of the above-mentioned equation were discussed through this method separately such as convection-diffusion-wave equation \eqref{4.1},  reaction-diffusion-wave equation \eqref{4.2} and diffusion-wave equation \eqref{4.3}. Moreover, we explained how to derive the exact solutions for the underlying equation along with initial conditions using the obtained invariant subspaces. Finally, we extended this method to two-dimensional time-fractional non-linear PDEs with time delay \eqref{D1}.  Also, the effectiveness and applicability of the method were illustrated through the two-dimensional time-fractional convection-reaction-diffusion-wave equation with time delay \eqref{6.9}. In addition, we observe that the obtained exact solutions can be viewed as the combinations of Mittag-Leffler function and polynomial, exponential and trigonometric type functions. We would like to point out that the obtained results are new and interesting, also observe that the given equation has not been discussed anywhere in the literature. This study shows that the discussed method and the obtained results are a very useful and efficient mathematical method to derive exact solutions for various types of integer-order and non-integer scalar and coupled non-linear PDEs in the fields of science and engineering for future research.

\end{document}